\pdfoutput=1
\documentclass[11pt]{amsart}

\usepackage[T1]{fontenc}
\usepackage[utf8]{inputenc}

\usepackage{amsmath,amssymb,amsthm,mathtools}
\usepackage{graphicx}
\usepackage{xcolor}
\usepackage[sort,compress]{cite}
\usepackage{tikz-cd}
\usepackage{changepage}
\usepackage[hidelinks]{hyperref}
\usepackage[
    a4paper,
    margin=2.5cm
]{geometry}

\hypersetup{
  pdftitle={Variational Lifting and Optimal Gauges on Riemannian Homogeneous Spaces},
  pdfauthor={J. R. Goodman and L. J. Colombo}
}

\theoremstyle{plain}
\newtheorem{theorem}{Theorem}[section]
\newtheorem{lemma}[theorem]{Lemma}
\newtheorem{corollary}[theorem]{Corollary}

\theoremstyle{definition}
\newtheorem{definition}[theorem]{Definition}

\theoremstyle{remark}
\newtheorem{remark}[theorem]{Remark}

\newcommand{\CP}{\mathbb{CP}}
\newcommand{\su}{\mathfrak{su}}
\newcommand{\Horsp}{\operatorname{Hor}}
\newcommand{\Versp}{\operatorname{Ver}}

\newcommand{\so}{\mathfrak{so}}

\newcommand{\g}{\mathfrak{g}}

\newcommand{\Exp}{\normalfont{\text{Exp}}}

\newcommand{\ext}{\normalfont{\text{ext}}}

\newcommand{\comment}[1]{}
\newcommand{\todo}[1]{}

\makeatletter
\def\widebreve{\mathpalette\wide@breve}
\def\wide@breve#1#2{\sbox\z@{$#1#2$}%
     \mathop{\vbox{\m@th\ialign{##\crcr
\kern0.08em\brevefill#1{0.8\wd\z@}\crcr\noalign{\nointerlineskip}%
                    $\hss#1#2\hss$\crcr}}}\limits}
\def\brevefill#1#2{$\m@th\sbox\tw@{$#1($}%
  \hss\resizebox{#2}{\wd\tw@}{\rotatebox[origin=c]{90}{\upshape(}}\hss$}
\makeatother

\newcommand{\Z}{\mathbb{Z}}

\newcommand{\R}{\mathbb{R}}
\newcommand{\grad}{\normalfont{\text{grad}}}

\newcommand{\tr}{\normalfont{\text{tr}}}
\newcommand{\SO}{\normalfont{\text{SO}}}

\newcommand{\ad}{\normalfont{\text{ad}}}
\newcommand{\Ad}{\normalfont{\text{Ad}}}

\begin{document}

\title{Variational Lifting and Optimal Gauges on Riemannian Homogeneous Spaces}

\author{Jacob R. Goodman}
\address{Norwegian University of Science and Technology, Norway}
\email{jacob.goodman@ntnu.no}

\author{Leonardo J. Colombo}
\address{Centre for Automation and Robotics (CSIC-UPM),\\
Ctra. M-300, Campo Real, Km 0.200, 28500 Arganda del Rey, Madrid, Spain}
\email{leonardo.colombo@car.upm-csic.es}

\subjclass[2020]{70H33, 70G75, 53C30}

\keywords{Riemannian homogeneous spaces; variational lifting; gauge invariance; Euler-Poincar\'e reduction; symmetry breaking}

\begin{abstract}
Building on standard Euler-Poincar\'e reduction on Lie groups, we develop a variational lifting framework for mechanical systems on Riemannian homogeneous spaces $H=G/K$. A mechanical action on $H$ is lifted to $G$ through a functional whose kinetic energy depends only on the horizontal component of the velocity, and we prove that its critical points project precisely onto those of the original action.

The lifted functional possesses a natural gauge invariance under $H^1$ curves with values in the isotropy subgroup $K$. Consequently, every lifted critical curve decomposes into a smooth horizontal representative and an arbitrary vertical gauge, and the associated Euler-Poincar\'e equations with symmetry breaking are obtained in reduced form. We then introduce a second variational problem that selects a distinguished lift by minimizing the vertical kinetic energy.

The optimal gauge is a length-minimizing geodesic on $K$, yielding an explicit decomposition of the total energy into projected and vertical contributions and a geometric interpretation in terms of holonomy for closed projected curves. The framework is illustrated first for pure quantum states on $\mathbb{CP}^{n}\cong SU(n+1)/S(U(1)\times U(n))$, where the gauge freedom corresponds to the isotropy of unitary lifts, and then for $S^2\cong SO(3)/SO(2)$ through an optimal camera-orientation problem for an axisymmetric satellite subject to an undesirable pointing region.
\end{abstract}

\maketitle

\section{Introduction}

Symmetry reduction is one of the organizing principles of geometric mechanics.
Its classical origins go back to Euler's equations for the free rigid body
\cite{Euler1765} and to Poincar\'e's formulation of mechanical equations in
variables adapted to a Lie-group action \cite{Poincare1901}. Arnold later
interpreted the Euler equations as geodesic equations for invariant metrics on
Lie groups and extended this viewpoint to ideal fluids \cite{Arnold1966}.
In modern form, Lagrangian reduction and the Euler--Poincar\'e equations
provide the variational counterpart of Lie--Poisson reduction; see, for
instance,
\cite{CendraHolmMarsdenRatiu1998,HolmMarsdenRatiu1998,MarsdenRatiu,Holm}.
In particular, parameter-dependent Lagrangians and advected quantities provide
a natural mechanism for retaining reduction when the full symmetry of the
mechanical system is broken \cite{HolmMarsdenRatiu1998}.

Full symmetry is frequently destroyed by the data defining a mechanical or
control problem. A preferred direction, an external field, a target
configuration, or an obstacle may break the invariance of a potential or cost
while preserving an isotropy subgroup. Such systems may nevertheless be
reduced by adjoining an advected parameter and recovering invariance on an
enlarged space. This mechanism appears in classical examples such as the heavy
top and has been developed for continuum and complex-fluid models
\cite{HolmMarsdenRatiu1998,GayBalmazTronci2010}, higher-order variational
problems \cite{GayBalmazHolmMeierRatiuVialard2012}, and optimal control
problems with symmetry-breaking costs
\cite{BlochColomboGuptaOhsawa2017}. For homogeneous spaces, Lagrangian
reduction with advected parameters was studied in \cite{Vizman2015}. A broad
recent account of broken symmetry, emphasizing the role of successive
compositions of geometric maps, is given by Holm in
\cite{HolmGeometricMechanicsIII}.

A complementary variational viewpoint is provided by Clebsch formulations.
Clebsch variables and momentum maps allow Euler--Poincar\'e dynamics to be
obtained from constrained variational principles on enlarged spaces; this
perspective is closely related to the role of semidirect products in geometric
mechanics \cite{MarsdenRatiuWeinstein1984}. In optimal control, Gay-Balmaz and
Ratiu \cite{GayBalmazRatiu2011} used the Clebsch approach to formulate a broad
class of problems including rigid-body dynamics, ideal flows, and
double-bracket systems. The construction developed here belongs to the same
variational tradition, but the enlargement arises from a different source:
the nonuniqueness of lifting a curve from a homogeneous space to its Lie group.

Let $H=G/K$ be a Riemannian homogeneous space with canonical projection
$\pi:G\rightarrow H$. A curve $q(t)\in H$ admits infinitely many lifts:
if $g(t)$ is one lift, then $g(t)k(t)$ has the same projection for every
$K$-valued curve $k(t)$. This is the natural gauge freedom associated with the
principal bundle $G\rightarrow G/K$. In the variational problem considered
here this freedom becomes an actual invariance of the lifted functional:
both the lifted potential and the horizontal kinetic energy are insensitive to
time-dependent vertical transformations. Hence the choice of lift is not
merely a choice of coordinates or representative; the variational principle
itself is degenerate along the isotropy directions.

This interpretation also fits naturally with the
composition-of-maps viewpoint emphasized in
\cite{HolmGeometricMechanicsIII}. Two distinct compositions appear in our
framework. First, the advected parameter
$\alpha=\Psi_{g^{-1}}(\alpha_0)$ describes a fixed reference datum as seen
through the inverse group motion and encodes the broken symmetry. Second,
every critical lift factors as $g=hk$, where the horizontal factor $h$
determines the projected motion and $k\in K$ modifies only its representative
in the fiber. These mechanisms play complementary roles: the former restores
the symmetry required for reduction, whereas the latter exposes the gauge
freedom created by lifting from $G/K$ to $G$. The additional variational
principle introduced below then selects a distinguished representative inside
this gauge class.

Variational reduction on Lie groups and homogeneous spaces was previously
used by the present authors in the obstacle-avoidance problem studied in
\cite{GoodmanColombo2024}. That work concerns a higher-order variational
problem leading to Riemannian cubic-type trajectories and derives reduced
necessary conditions on homogeneous spaces through a connection on the
horizontal bundle. The question addressed here is different. We consider the
ordinary mechanical action and use the homogeneous-space structure to study
the variational problem satisfied by its lifts. The Euler--Poincar\'e material
on Lie groups is recalled only to establish the reduction with partial
symmetry required later. The new ingredients are the variational lifting
itself, its gauge invariance, the regularity consequences of this degeneracy,
and a second minimum-energy variational problem selecting an optimal gauge.
In particular, the $H^1$ setting is intrinsic to the construction: arbitrary
$H^1$ curves in $K$ are genuine gauges of the lifted functional, so a lifted
critical curve need not be smooth even though its projection and horizontal
representative are smooth.

More precisely, we prove that a mechanical action on $H$ can be lifted to a
functional on $G$ whose kinetic term depends only on the horizontal component
of the velocity, and that criticality of the two variational problems is
equivalent under projection. When the lifted potential admits a partial
symmetry, the corresponding dynamics satisfy reduced Euler--Poincar\'e
equations with an advected parameter. Every lifted critical point then
decomposes into a smooth horizontal representative and an arbitrary $H^1$
vertical gauge. We subsequently consider all lifts of a fixed projected
critical curve satisfying prescribed group endpoints and minimize their full
kinetic energy. This problem reduces exactly to a geodesic problem on the
isotropy subgroup $K$. Consequently, the optimal gauge is a
length-minimizing geodesic and the total energy splits into the projected
energy and an explicit minimum vertical contribution. For closed projected
curves, this contribution is the minimum energy required to compensate the
holonomy of the horizontal lift.

Two examples illustrate different aspects of the framework. For pure quantum
states,
$\CP^n\cong SU(n+1)/S(U(1)\times U(n))$, unitary evolutions are lifts of
curves in projective Hilbert space, while the isotropy subgroup represents
unitary motions that leave the physical state unchanged. A fidelity potential
provides a natural symmetry-breaking term, and the optimal gauge determines
the minimum additional unitary motion compatible with prescribed unitary
endpoints without altering the physical state trajectory. For
$S^2\cong SO(3)/SO(2)$, we consider a camera mounted on an axisymmetric
satellite. A symmetry-breaking potential penalizes an undesirable pointing
region, the projected variational problem determines the camera trajectory on
$S^2$, and the remaining $SO(2)$ freedom is resolved by the minimum-energy
gauge. The numerical illustration uses a Lie group variational integrator on
$SO(3)$ to display both the projected motion and the resulting optimal lift.

The paper is organized as follows. Section~2 recalls the Euler--Poincar\'e
reduction with partial symmetry needed for the subsequent construction.
Section~3 develops the variational lifting on Riemannian homogeneous spaces,
establishes its gauge invariance and regularity properties, and solves the
minimum-energy gauge problem; the quantum example concludes that section.
Section~4 applies the framework to the satellite camera-pointing problem and
presents the numerical illustration. Section~5 contains conclusions and
directions for future work.

\section{Euler-Poincar\'e Reduction with Broken Symmetry}\label{Section: Euler-Poincare}\label{subsection: EP broken symmetry}

Let $G$ be a Lie group with identity element $e\in G$ and Lie algebra $\g = T_e G$. Assume that $G$ is equipped with a left-invariant Riemannian metric $\langle \cdot, \cdot \rangle$ and Levi-Civita connection $\nabla$. Denote by $H^1([a,b],G)$ the Sobolev space of absolutely
continuous curves whose weak derivative is square integrable, understood
in local coordinates. It is well-known that the Lie algebra $\g$ is canonically isomorphic to the space of left-invariant vector fields $\mathfrak{X}_L(G)$ on $G$ via the map:
$$\phi: \g \to \mathfrak{X}_L(G), \qquad \phi(\xi)(g) = T_eL_g(\xi)$$
where $L_g: G \to G$ denotes left-translation by $g \in G$. This isomorphism allows to construct an operator $\nabla^{\g}: \g \times \g \to \g$ defined by:
\begin{equation}\label{eq: g-connection}
    \nabla^\g_\xi \eta := \nabla_{\phi(\xi)}\phi(\eta)(e),
\end{equation}
for all $\xi, \eta \in \g$. We refer to $\nabla^{\mathfrak g}$ as the Riemannian
$\mathfrak g$-connection associated with $\nabla$. 

Equivalently, we have the following representation that we will make extensive use of (see Theorem~5.40 of \cite{bullo2019geometric} for more details):
\begin{lemma}\label{lemma: covg-decomp}
For all $\xi, \eta \in \g$ the Riemannian $\g$-connection can be expressed as \begin{align*}
    \nabla_{\xi}^\g \eta = \frac12\left( [\xi, \eta] - \ad^\dagger_{\xi} \eta - \ad^\dagger_{\eta} \xi\right).
    \end{align*}
Here
$[\cdot,\cdot]$ denotes the Lie bracket on $\mathfrak g$, and
$\ad_\xi:\mathfrak g\to\mathfrak g$ is defined by
$\ad_\xi\eta=[\xi,\eta]$. We also denote by $\ad_\xi^\dagger$ its metric
adjoint, characterized by $\langle \ad_\xi^\dagger\eta,\sigma\rangle= \langle\eta,[\xi,\sigma]\rangle,
    \, \xi, \eta,\sigma\in\mathfrak g$.
\end{lemma}

% Although $\nabla^{\g}$ is not a connection, we shall refer to it as the \textit{Riemannian $\g$-connection} corresponding to $\nabla$ because of the similar properties that it satisfies (see \cite{goodman2022reduction}):

% We define an operator $\ad^\dagger: \g \times \g \to \g$ by $\ad^\dagger_{\xi} \eta = (\ad^\ast_{\xi} \eta^\flat)^\sharp$ for all $\xi, \eta \in \g$. It follows that $\ad^\ast$ is bilinear, and satisfies
% $\left<\ad^\dagger_\xi \eta, \sigma \right> = \left<\ad^\ast_\xi \eta^\flat, \sigma \right>_\ast = \left<\eta^\flat, \ad_\xi \sigma \right>_\ast = \left<\eta, \ad_\xi \sigma \right>$ for all $\xi, \eta, \sigma \in \g$. Hence, $\ad^\ast$ is nothing more than the adjoint of the adjoint action with respect to the Riemannian metric on $G$. This leads to the following decomposition of $\nabla^\g$ (see Theorem $5.40$ of \cite{bullo2019geometric}, for instance)

For a $C^1$ curve $g:[a,b]\to G$ and a $C^1$ vector field $X$ along $g$, we denote by
$D_tX:=\nabla_{\dot g}X$ the covariant derivative of $X$ along $g$.

For $X_g\in T_gG$, we use the shorthand
$g^{-1}X_g:=T_gL_{g^{-1}}(X_g)$ for its left trivialization.
The operator $\nabla^\g$ can be understood as measuring the failure
of left-translation to commute with covariant differentiation.
The following standard moving-frame identity, whose smooth version
appears in~\cite{GoodmanThesis}, describes this relation. We also
record its weak $H^1$ extension, which will be needed below.

\begin{lemma}\label{lemma: cov-to-covg}
Let $g:[a,b]\to G$ be a $C^1$ curve and $X$ be a $C^1$ vector field
along $g$. Then the following relation holds for all $t\in[a,b]$:
\[
D_tX(t)
=
T_eL_{g(t)}
\left(
\dot\eta(t)+\nabla^{\mathfrak g}_{\xi(t)}\eta(t)
\right),
\]
where
\[
\xi(t):=
T_{g(t)}L_{g(t)^{-1}}\bigl(\dot g(t)\bigr),
\qquad
\eta(t):=
T_{g(t)}L_{g(t)^{-1}}\bigl(X(t)\bigr).
\]
The same identity holds in the weak sense for $H^1$ curves $g$ and
$H^1$ vector fields $X$ along $g$.
\end{lemma}

\begin{proof}
The $C^1$ statement is proved in~\cite{GoodmanThesis}. We prove here
the weak $H^1$ extension. 

Choose smooth approximations
$g_n\to g$ and $\eta_n\to\eta$ in $H^1$, and set
$X_n:=T_eL_{g_n}(\eta_n)$. The $C^1$ identity holds for each
$(g_n,X_n)$.  

Since $H^1([a,b])\hookrightarrow C^0([a,b])$,
the convergence $g_n\to g$ and $\eta_n\to\eta$ is uniform, while
$\dot g_n\to\dot g$ and $\dot\eta_n\to\dot\eta$ in $L^2$. Since the group operations and left translations are smooth and
$\nabla^{\mathfrak g}$ is bilinear, all terms converge in $L^2$.
Passing to the limit gives the same identity in the weak sense.
\end{proof}

 Given $g_a,g_b\in G$ and $a<b$, let $\Omega^{a,b}_{g_a,g_b}(G)$ denote the space of $H^1$ curves $g:[a,b]\to G$ satisfying $g(a)=g_a$ and $g(b)=g_b$. Its tangent space consists of $H^1$ vector fields along $g$ vanishing at the endpoints. We will frequently use the simplified notation $\Omega$ when the endpoints, time interval, and space are clear from context. We refer to $\Omega$ as the admissible path space, and elements $g \in \Omega$ as admissible curves. The $H^1$ setting will be essential for the lifted variational problem on homogeneous spaces, where arbitrary $H^1$ curves with values in the isotropy subgroup arise naturally as vertical gauges. %A curve $g_s(s, \cdot) \in \Omega$ with $s \in (-\varepsilon, \varepsilon)$ for some $\varepsilon > 0$ and $g := g_s(0, \cdot)$ is called an admissible variation along $g$, and the vector field $\delta g = \partial_s g_s(0, \cdot)$ is called the variational vector field of $g_s$. See \cite{Hebey} for more information on Sobolev spaces on Riemannian manifolds.
We regard $\Omega$ as the corresponding Hilbert submanifold of
$H^1([a,b],G)$. 

An admissible variation of $g\in\Omega$ is a $C^1$ curve
\[
(-\varepsilon,\varepsilon)\ni s\longmapsto g_s\in\Omega,
\qquad g_0=g.
\]
Its variational vector field is
\[
\delta g
:=
\left.\frac{d}{ds}\right|_{s=0}g_s
\in T_g\Omega.
\]
Time derivatives of $H^1$ curves, and the identities involving them below, are understood in the weak sense and almost everywhere until additional regularity has been established. In local coordinates, this implies that the variational vector
field belongs to $H^1$ and that mixed derivatives are understood in the
distributional sense.

Consider the action functional $\mathcal A:\Omega\to\mathbb R$ given by
\begin{equation}\label{eq: mechanical lagrangian}
\mathcal A[g]
=
\int_a^b
\left(
\frac12\|\dot g(t)\|^2-V(g(t))
\right)dt,
\end{equation}
where $\|\cdot\|$ denotes the norm induced by the left-invariant
Riemannian metric on $G$, and $V:G\to\mathbb R$ is a smooth scalar
field called the potential. This is the usual action associated to the
mechanical Lagrangian $L(g,\dot g)=\frac12\|\dot g\|^2-V(g)$.  We note that the potential may be a potential energy corresponding to a physical system, or an \textit{artificial potential} used to modify the behavior of the resulting minimizers, such as in the case of the obstacle avoidance problem (see \cite{BlCaCoCDC}). The critical points of Equation~\eqref{eq: mechanical lagrangian} over the admissible path space satisfy the Euler-Lagrange equations in the weak sense. For a mechanical Lagrangian, these equations take the form
\begin{equation}\label{eq: EL equations mechanical lagrangian}
D_t \dot{g}(t) = -\grad V(g(t)).
\end{equation} where $\grad V$ denotes the Riemannian gradient of $V$ with respect to
the left-invariant metric on $G$. This identity holds almost everywhere on $[a,b]$; as shown below, critical points are smooth, and the equation therefore holds pointwise after the corresponding regularity bootstrap.

Writing $\xi(t):=T_{g(t)}L_{g(t)^{-1}}\big(\dot g(t)\big)$, the left-invariance of the metric identifies the kinetic part of the dynamics with an equation on the Lie algebra $\mathfrak g$, together with the reconstruction equation $\dot g(t)=T_eL_{g(t)}\big(\xi(t)\big)$. In the absence of a potential, or when the potential is invariant under left translations, this produces a closed variational principle in the reduced variable $\xi$. For a general potential, however, the
left-trivialized force
$T_{g(t)}L_{g(t)^{-1}}\bigl(\grad V(g(t))\bigr)$
still depends on the group variable $g$, and therefore the equations
do not necessarily close in the reduced variable $\xi$ alone. The reduction procedure can nevertheless be recovered when this dependence is encoded by an advected parameter. This is the situation described in the following definition.

\begin{definition}\label{def: partial symmetry Q}
We say that the potential $V: G \to \R$ admits a \textit{partial symmetry} over some smooth manifold $M$ if there exists a smooth left-action $\Psi: G \times M \to M$ and a smooth function $V_{\ext}: G \times M \to \R$ called the \textit{extended potential} such that:
\begin{enumerate}
    \item There exists a fixed $\alpha_0\in M$ such that, for every $g\in G$, we have $V_{\ext}(g, \alpha_0) = V(g)$.
    \item $V_{\ext}$ is invariant under the left-action $\chi: G \times (G \times M) \to G \times M$ given by $\chi_h(g, \alpha) = (hg, \Psi_{h}(\alpha))$, that is, $V_{\ext}(hg,\Psi_h(\alpha))=V_{\ext}(g,\alpha)$ for all $g, h \in G, \alpha \in M$.
\end{enumerate}
\end{definition}

If $G_{\alpha_0}:=\{h\in G:\Psi_h(\alpha_0)=\alpha_0\}$ denotes the isotropy subgroup of $\alpha_0$, then the above conditions
imply
\[
V(hg)=V(g),
\qquad h\in G_{\alpha_0},\quad g\in G.
\]
Thus $V$ retains the symmetry associated with the isotropy subgroup
$G_{\alpha_0}$, while the full $G$-symmetry is recovered by the
extended potential.

The left-action $\Psi$ induces a natural infinitesimal action of the Lie algebra on $M$. 
For any $\xi \in \g$, the \textit{infinitesimal generator} associated to $\xi$ is the vector field $\xi_M \in \mathfrak{X}(M)$ defined by
\begin{equation}\label{eq: infinitesimal_gen}
    \xi_M(x) = \frac{d}{dt} \Big{\vert}_{t=0} \Psi_{\Exp(t\xi)}(x)
,
\end{equation}
for all $x \in M$, where $\Exp: \g \to G$ is the Lie exponential map on $G$.

\begin{theorem}\label{prop: EP geo symmetry breaking}
Let $G$ be a Lie group equipped with a left-invariant metric, and suppose that $g: [a, b] \to G$ is a critical point of the action functional defined in \eqref{eq: mechanical lagrangian}. Then $g$ is smooth, and if the potential admits a partial symmetry over some smooth manifold $M$, then $\xi(t):=T_{g(t)}L_{g(t)^{-1}}\big(\dot g(t)\big)$ and $\alpha(t):=\Psi_{g(t)^{-1}}(\alpha_0)$ satisfy
\begin{align}
    \dot{\xi} + \nabla^\g_{\xi} \xi &= -\grad_1 V_{\ext}(e, \alpha),\label{EP: mechanical geo} \\
    \dot{\alpha} &= -\xi_M(\alpha), \label{EP: parameter}
\end{align} where $\grad_1V_{\rm ext}(g,\alpha)$ denotes the Riemannian
gradient of $V_{\rm ext}$ with respect to its first argument.
\end{theorem}

 \begin{proof}
Suppose that $g_s$ is an admissible variation of $g$ with left-trivialized variational vector field $\eta(t):=
T_{g(t)}L_{g(t)^{-1}}\big(\delta g(t)\big)$.
Conversely, every
$\eta\in H^1_0([a,b],\mathfrak g)$ is generated by the admissible
variation $g_s(t)=g(t)\Exp(s\eta(t))$, for $|s|$ sufficiently small.

If $g$ is a critical point of the action \eqref{eq: mechanical lagrangian}, then:
\begin{align*}
    0 = \frac{\partial}{\partial s}\Big{\vert}_{s=0} \int_a^b \left(\frac12 \|\dot{g}_s\|^2 - V(g_s) \right)dt = \int_a^b \left( \left<D_s \dot{g}_s, \dot{g}_s\right> - \frac{\partial}{\partial s}V(g_s)\right)\Big{\vert}_{s=0} dt.
\end{align*} Since the Levi--Civita connection is torsion-free and the parameter vector
fields of the variation commute,
$\left.D_s\dot g_s\right|_{s=0}=D_t\delta g$. Using Lemmas~\ref{lemma: cov-to-covg} and~\ref{lemma: covg-decomp}, and interpreting all derivatives in the
weak sense, we obtain
\[
\int_a^b
    \langle D_t\delta g,\dot g\rangle\,dt
=
\int_a^b
\left(
    \langle \xi,\dot\eta\rangle
    +
    \langle \ad^\dagger_\xi\xi,\eta\rangle
\right)dt.
\]
For a general smooth potential, define
\[
    F(g):=T_gL_{g^{-1}}\bigl(\grad V(g)\bigr).
\]
Then
\[
    \delta\mathcal A[g]
    =
    \int_a^b
    \left(
        \langle \xi,\dot\eta\rangle
        +
        \langle \ad^\dagger_\xi\xi-F(g),\eta\rangle
    \right)dt
    =0
\]
for every $\eta\in H^1_0([a,b],\mathfrak g)$.
Hence, by the definition of the distributional derivative,
\[
    \dot\xi
    =
    \ad^\dagger_\xi\xi-F(g)
\]
in the sense of distributions.

Since $g\in H^1([a,b],G)$, its left-trivialized velocity satisfies
$\xi\in L^2([a,b],\g)$. Moreover, in one dimension
$H^1([a,b],G)\subset C^0([a,b],G)$, and since $V$ is smooth, the map
$g\mapsto T_gL_{g^{-1}}\big(\grad V(g)\big)$
is smooth. Hence $F(g)\in C^0([a,b],\g)$. Since $\ad^\dagger$ is a bounded bilinear map on the finite-dimensional Lie algebra $\g$, we have
$\ad^\dagger_\xi\xi\in L^1([a,b],\g)$. The weak equation therefore implies
$\dot\xi\in L^1([a,b],\g)$,
and consequently
$\xi\in W^{1,1}([a,b],\g)\subset C^0([a,b],\g)$.
The reconstruction equation
$\dot g(t)=T_eL_{g(t)}\big(\xi(t)\big)$
then implies that $g\in C^1([a,b],G)$. It follows that both $F(g)$ and $\ad^\dagger_\xi\xi$ are continuous,
and therefore the weak equation gives
$\xi\in C^1([a,b],\mathfrak g)$. The reconstruction equation then
gives $g\in C^2([a,b],G)$. Repeating this argument and using the
smoothness of the coefficients yields, by bootstrap, that both
$\xi$ and $g$ are smooth.

Moreover, if $V$ admits a partial symmetry, then
\begin{align*}
\int_a^b \frac{\partial}{\partial s}V(g_s)\Big{\vert}_{s=0} dt
&= \int_a^b \frac{\partial}{\partial s}V_{\ext}(g_s,\alpha_0)\Big{\vert}_{s=0} dt \\
&= \int_a^b \frac{\partial}{\partial s}V_{\ext}(g^{-1}g_s,\alpha)\Big{\vert}_{s=0} dt \\
&= \int_a^b \left<\grad_1 V_{\ext}(e,\alpha),\eta\right>dt.
\end{align*}
Since $\eta$ is arbitrary, comparing this expression with the definition of $F(g)$ and applying the fundamental lemma of the calculus of variations gives
$F(g)=\grad_1V_{\ext}(e,\alpha)$.
 Therefore,
\[
\dot{\xi}
=
\ad^\dagger_\xi\xi-\grad_1V_{\ext}(e,\alpha),
\]
which is equivalent to Equation~\eqref{EP: mechanical geo} by Lemma~\ref{lemma: covg-decomp}.

Finally, differentiating
$\alpha(t)=\Psi_{g(t)^{-1}}(\alpha_0)$ and using
$g(t+\varepsilon)^{-1}
    =
    \Exp(-\varepsilon\xi(t))g(t)^{-1}
    +o(\varepsilon)$, together with the definition of the infinitesimal generator, gives
$\dot\alpha=-\xi_M(\alpha)$, which proves Equation~\eqref{EP: parameter}.
\end{proof}

\section{Variational Lifting on Riemannian Homogeneous Spaces}\label{sec6}

\subsection{Variational lifting to the Lie group}
We recall only the notation needed for the reduction theorem; additional background on homogeneous spaces can be found in \cite{helgason2024geometric}. Let $G$ be a connected Lie group, $K\subset G$ be a closed Lie subgroup,
and $H:=G/K$ with canonical projection $\pi:G\to H$. Let
$\mathfrak{s}:=T_eK$ and assume that $G$ is equipped with a
left-invariant Riemannian metric whose inner product on $\mathfrak{g}$
is $\operatorname{Ad}_K$-invariant, namely,
\[
\langle \operatorname{Ad}_k\xi,\operatorname{Ad}_k\eta\rangle_G
=
\langle \xi,\eta\rangle_G,
\qquad
k\in K,\quad \xi,\eta\in\mathfrak{g}.
\]
Equivalently, the left-invariant Riemannian metric on $G$ is also
right-$K$-invariant. Let
$\mathfrak{h}:=\mathfrak{s}^{\perp}$. Since $\mathfrak{s}$ is
$\operatorname{Ad}_K$-invariant, the $\operatorname{Ad}_K$-invariance
of the metric implies that $\mathfrak{h}$ is also
$\operatorname{Ad}_K$-invariant. We equip $H$ with the induced
$G$-invariant Riemannian metric, so that
$\pi:G\to H$ is a Riemannian submersion. Explicitly, this metric is characterized by
\[
\left\langle
T_g\pi\bigl(T_eL_g\xi\bigr),
T_g\pi\bigl(T_eL_g\eta\bigr)
\right\rangle_H
=
\langle\xi,\eta\rangle_G,
\hbox{ for all } \xi,\eta\in\mathfrak h \hbox{ and }g\in G
\]
which is well defined by the $\Ad_K$-invariance assumptions.

Since $\pi$ is a submersion, we obtain a canonical vertical space
\[
\Versp_g
:=
\ker(T_g\pi)=T_eL_g(\mathfrak{s}).
\]
Taking the disjoint union of all vertical spaces, we obtain the
vertical bundle
\[
VG:=\bigsqcup_{g\in G}\Versp_g.
\]

The horizontal space is the orthogonal complement of the vertical
space and, by left-invariance of the metric, satisfies

\[
\Horsp_g:=\Versp_g^\perp=T_eL_g(\mathfrak h)
\]
for all $g\in G$. The horizontal bundle is
\[
HG:=\bigsqcup_{g\in G}\Horsp_g.
\]
In such a case, the tangent bundle of $G$ decomposes as
\begin{equation}\label{eq:direct_sum_TG}
TG=VG\oplus HG.
\end{equation}

The right-$K$-invariance of the metric, equivalently the
$\operatorname{Ad}_K$-invariance above, implies that the horizontal
distribution is right-$K$-equivariant,
\[
T_gR_k(\operatorname{Hor}_g)=\operatorname{Hor}_{gk},
\qquad
g\in G,\quad k\in K,
\]
and therefore defines a principal connection on $G\to G/K$.

We denote by $\mathcal V,\mathcal H:TG\to TG$ the vertical and
horizontal projections associated with
\eqref{eq:direct_sum_TG}, respectively; thus, for
$X_g\in T_gG$,
\[
X_g=\mathcal V(X_g)+\mathcal H(X_g),
\qquad
\mathcal V(X_g)\in\Versp_g,
\qquad
\mathcal H(X_g)\in\Horsp_g.
\]

Considering~\eqref{eq:direct_sum_TG} at the identity $e\in G$, the
Lie algebra decomposes as
$\mathfrak g=\mathfrak s\oplus\mathfrak h$, where
$\mathfrak s$ is the Lie algebra of $K$.

The restriction
$T_e\pi|_{\mathfrak h}:\mathfrak h\to T_{\pi(e)}H$
is an isomorphism. Thus, under the chosen orthogonal decomposition,
\[
\mathfrak h\cong\mathfrak g/\mathfrak s
\cong T_{\pi(e)}H.
\]

For each $g\in G$, we denote the restrictions of the projections by
\[
\mathcal V_g
:=
\mathcal V|_{T_gG}
:
T_gG\to\Versp_g,
\qquad
\mathcal H_g
:=
\mathcal H|_{T_gG}
:
T_gG\to\Horsp_g.
\]
%A section of the vertical (horizontal) bundle is called a vertical (horizontal) vector field. \textcolor{blue}{A vector field $Y\in\mathfrak X(G)$ is called projectable if it is $\pi$-related to a vector field $X\in\mathfrak X(H)$, that is,
%\[
%T_g\pi\big(Y(g)\big)=X\big(\pi(g)\big)
%\]
%for every $g\in G$. A basic vector field is a horizontal projectable vector field.} The space of basic vector fields is denoted by $\mathcal{B}(G)$. Given a vector field $X \in \mathfrak{X}(H)$, we denote its horizontal lift by $\tilde{X} \in \mathcal{B}(G)$.

%We assume throughout this section that $H$ is a Riemannian homogeneous space of $G$, so that $\pi:G\to H$ is a Riemannian submersion. Denote the Levi-Civita connections on $H$ and $G$ by $\nabla$ and $\tilde{\nabla}$, respectively. The Riemannian metric on $G$ is assumed to be left-invariant unless otherwise stated.

%\begin{lemma}\label{lemma: cov_G_to_H}
%Let $X, Y \in \Gamma(TH)$, and $\tilde{X}, \tilde{Y} \in \mathcal{B}(G)$ be the horizontal lifts of $X$ and $Y$, respectively. Further suppose that $V \in \Gamma(VG)$. Then, the following identities hold:
%\begin{align}
    %\tilde{\nabla}_{\tilde{X}}\tilde{Y}
    %&= \widetilde{\nabla_X Y}
    %+ \frac12\mathcal{V}([\tilde{X}, %\tilde{Y}]),
    %\label{eq: covG-to-covH}\\
    %\textcolor{blue}{
    %\mathcal H\big(\tilde{\nabla}_{\tilde X}V\big)
    %&=
    %\mathcal H\big(\tilde{\nabla}_{V}\tilde X\big)}
    %\label{cov_G_commute}
%\end{align}
%\end{lemma}

%Next, we seek a variational principle on $G$ which is equivalent to the variational principle $\delta\mathcal A_H=0$ on $H$. 

Let
\[
    \Omega(G):=\Omega^{a,b}_{g_a,g_b}(G),
    \qquad
    \Omega(H):=\Omega^{a,b}_{q_a,q_b}(H),
\]
where $\pi(g_a)=q_a$ and $\pi(g_b)=q_b$. We shall use the standard fact that, if
$g\in H^1([a,b],G)$ and $\pi:G\to H$ is smooth, then $q=\pi\circ g$ belongs to
$H^1([a,b],H)$ and
\[ \dot q(t)=T_{g(t)}\pi\big(\dot g(t)\big),
    \qquad
   \text{for a.e. }t\in[a,b].
\]
We denote by $\|\cdot\|_G$ and $\|\cdot\|_H$ the pointwise norms
induced by the Riemannian metrics on $G$ and $H$, respectively. Since $\pi:G\to H$ is a Riemannian submersion, it follows that
\[
    \|\dot q(t)\|_H
    =
    \big\|T_{g(t)}\pi\big(\dot g(t)\big)\big\|_H
    =
    \|\mathcal H_{g(t)}(\dot g(t))\|_G,
    \qquad
    \text{for a.e. }t\in[a,b].
\]%Throughout the section, \textcolor{blue}{we work under the quotient-metric assumptions stated above}. For some smooth scalar field $V: H \to \R$, we consider the action \textcolor{blue}{functional} $\mathcal{A}_H: \Omega(H) \to \R$ on $H$ as in \eqref{eq: mechanical lagrangian}. As usual, our aim is to find the critical points of $\mathcal{A}_H$ over $\Omega$. Since $H$ is a Riemannian manifold, this may be done via variational analysis. However, in many cases, it will be more convenient to work on $G$--especially in cases where we can apply reduction by symmetry.
 For a smooth potential $V:H\to\mathbb R$, consider the mechanical
action
\[
\mathcal A_H[q]
=
\int_a^b
\left(
\frac12\|\dot q(t)\|_H^2-V(q(t))
\right)dt,
\qquad q\in\Omega(H).
\]
We seek a variational principle on $G$ whose criticality condition
is equivalent to that of $\mathcal A_H$ along curves related by
projection.

This motivates the lifted action functional
\begin{equation}\label{eq: action function G geodesics}
    \tilde{\mathcal{A}}_G[g]
    =
    \int_a^b
    \left(
        \frac12 \|\mathcal{H}_{g(t)}(\dot{g}(t))\|^2_G
        -
        \tilde{V}(g(t))
    \right)dt,
\end{equation}
where $\widetilde V=V\circ\pi$.
We have the following Theorem:

\begin{theorem}\label{prop: variational principle equivalence homogeneous space geodesics}
Let $g\in\Omega(G)$ and set $q:=\pi\circ g\in\Omega(H)$. Then $q$ is a critical point of $\mathcal{A}_H$ if and only if $g$ is a critical point of $\tilde{\mathcal{A}}_G$.
\end{theorem}

\begin{proof}
For every $g\in\Omega(G)$, the definition of the lifted functional and the Riemannian-submersion property give
$\widetilde{\mathcal A}_G[g]
=
\mathcal A_H[\pi\circ g]$.
Suppose first that $q=\pi\circ g$ is critical for $\mathcal A_H$. If
$\delta g\in T_g\Omega(G)$, then
$\delta q(t)
:=
T_{g(t)}\pi\big(\delta g(t)\big)$ belongs to $T_q\Omega(H)$. Hence,
\[
\delta\widetilde{\mathcal A}_G[g](\delta g)
=
\delta\mathcal A_H[q](\delta q)
=
0,
\]
and therefore $g$ is a critical point of $\widetilde{\mathcal A}_G$.
Conversely, suppose that $g$ is critical for $\widetilde{\mathcal A}_G$, and let
$\delta q\in T_q\Omega(H)$ be arbitrary. Since
$\pi:G\to H$ is a Riemannian submersion, $\delta q$ has a unique horizontal
lift along $g$, given by
\[
\delta g(t)
:=
\left(
T_{g(t)}\pi\big|_{\Horsp_{g(t)}}
\right)^{-1}
\big(\delta q(t)\big).
\]
Since the horizontal-lift map depends smoothly on $g$, the $H^1$ regularity of $g$ and $\delta q$
implies that $\delta g\in H^1$. Moreover, since
$\delta q(a)=\delta q(b)=0$, we also have
$\delta g(a)=\delta g(b)=0$. Hence
$\delta g\in T_g\Omega(G)$. Therefore,
\[
\delta\mathcal A_H[q](\delta q)
=
\delta\widetilde{\mathcal A}_G[g](\delta g)
=
0.
\]
Since $\delta q$ is arbitrary, $q$ is a critical point of $\mathcal A_H$.
\end{proof}

From Theorem \ref{prop: variational principle equivalence homogeneous space geodesics}, we can generate critical points of $\mathcal{A}_H$ by finding critical points of $\tilde{\mathcal{A}}_G$ and then projecting them onto $H$. Notably, $g$ need not be a horizontal curve, since the kinetic energy in the lifted functional depends only on the horizontal component of $\dot g$. By left-trivializing the action \eqref{eq: action function G geodesics}, we obtain a reduced variational principle under the assumption that $\widetilde V$ admits a partial symmetry over some smooth manifold $M$, and subsequently obtain a form of Euler-Poincar\'e equations, provided that regularity can be established.

\subsection{Reduced Euler-Poincar\'e equations}

The critical points of $\mathcal A_H$ are smooth by the standard
one-dimensional regularity argument underlying Theorem \ref{prop: EP geo symmetry breaking}. However, the critical points of $\tilde{\mathcal{A}}_G$ are not generally smooth. Indeed, let $g\in\Omega(G)$ be a lift of a critical point $q$ of $\mathcal A_H$, and let $h$ be the horizontal lift of $q$ satisfying $h(a)=g(a)$. Then $h$ is smooth and, since $\pi\circ h=\pi\circ g$, $k(t):=h(t)^{-1}g(t)$ belongs to $K$, so that $k\in H^1([a,b],K)$. If $\widehat{k}\in H^1([a,b],K)$ is any other curve satisfying $\widehat{k}(a)=k(a)$ and $\widehat{k}(b)=k(b)$, then $\widehat{g}:=h\widehat{k}$ belongs to $\Omega(G)$ and satisfies $\pi(\widehat{g}(t))=\pi(h(t))=q(t)$.
By a slight abuse of notation, we also denote by
$\mathcal H:\mathfrak g\to\mathfrak h$ the orthogonal projection associated with $\mathfrak g=\mathfrak s\oplus\mathfrak h$. By left-invariance,
$\mathcal H_g\bigl(T_eL_g\xi\bigr)
=
T_eL_g\bigl(\mathcal H(\xi)\bigr)$. By Theorem \ref{prop: variational principle equivalence homogeneous space geodesics},
$\widehat{g}$ is therefore a critical point of $\tilde{\mathcal A}_G$. Since
$\widehat{k}$ may have only $H^1$ regularity, the lifted critical point
$\widehat{g}$ need not be smooth.

% {\color{red} \begin{definition}\label{def: basic variation}
% Let $g \in \Omega(G)$ be a basic curve. We call an admissible variation of $g$ through basic curves with a horizontal variational vector field a \textit{basic variation} of $g$. 
% \end{definition}
% }

\begin{theorem}\label{thm: EP_eqs_homo_G}
Suppose that $\widetilde V$ admits a partial symmetry over some smooth manifold $M$. Let $g\in\Omega(G)$, and define
\[
\xi:=g^{-1}\dot g,
\qquad
\alpha:=\Psi_{g^{-1}}(\alpha_0).
\]
Then $g$ is a critical point of $\widetilde{\mathcal A}_G$ if and only if the pair $(\xi,\alpha)$ satisfies the reduced variational principle
\[
\delta
\int_a^b
\left(
\frac12\left\|\mathcal H(\xi)\right\|_G^2
-
\widetilde V_{\ext}(e,\alpha)
\right)dt=0
\]
with respect to variations generated by
$\eta\in H^1([a,b],\mathfrak g)$
\[
\delta\xi=\dot\eta+[\xi,\eta],
\qquad
\delta\alpha=-\eta_M(\alpha),
\]
satisfying $\eta(a)=\eta(b)=0$.
In such a case, $\mathcal H(\xi)\in H^1([a,b],\mathfrak h)$, and $\xi$ and $\alpha$ satisfy
\begin{align}
\frac{d}{dt}\mathcal H(\xi)
&=
\ad^\dagger_\xi \mathcal H(\xi)-\grad_1\widetilde V_{\ext}(e,\alpha),
\label{eq: EP eqs horizontal geodesics}\\
\dot{\alpha}
&=
-\xi_M(\alpha).
\label{eq: parameter dynamics horizontal geodesics}
\end{align}
in the weak sense.
\end{theorem}

\begin{proof}
Since the horizontal distribution and the metric are left-invariant,
\[
\mathcal H_{g(t)}(\dot g(t))
=
T_eL_{g(t)}\bigl(\mathcal H(\xi(t))\bigr),
\qquad
\xi=g^{-1}\dot g,
\]
and therefore
\[
\left\|\mathcal H_{g(t)}(\dot g(t))\right\|_G
=
\left\|\mathcal H(\xi(t))\right\|_G
\]
for almost every $t\in[a,b]$. Moreover, by partial symmetry,
\[
\widetilde V(g)
=
\widetilde V_{\ext}(g,\alpha_0)
=
\widetilde V_{\ext}(e,\Psi_{g^{-1}}(\alpha_0))
=
\widetilde V_{\ext}(e,\alpha).
\]
Thus the lifted action becomes
\[
\widetilde{\mathcal A}_G[g]
=
\int_a^b
\left(
\frac12\left\|\mathcal H(\xi)\right\|_G^2
-
\widetilde V_{\ext}(e,\alpha)
\right)dt.
\]

Let $g_s$ be an admissible variation of $g$, and set
$\eta:=g^{-1}\delta g$. Then the standard left-trivialized variation formulas give
\[
\delta\xi=\dot\eta+[\xi,\eta],
\qquad
\delta\alpha=-\eta_M(\alpha),
\qquad
\eta(a)=\eta(b)=0.
\]
Conversely, every $\eta\in H^1([a,b],\mathfrak g)$ with
$\eta(a)=\eta(b)=0$ is generated by the admissible variation
$g_s(t)=g(t)\Exp(s\eta(t))$ for $|s|$ sufficiently small.
Differentiating the invariance identity
\[
\widetilde V_{\ext}
\bigl(\Exp(s\eta),\Psi_{\Exp(s\eta)}(\alpha)\bigr)
=
\widetilde V_{\ext}(e,\alpha)
\]
at $s=0$ gives
\[
d_2\widetilde V_{\ext}(e,\alpha)[\eta_M(\alpha)]
=
-d_1\widetilde V_{\ext}(e,\alpha)[\eta].
\]
Since $\delta\alpha=-\eta_M(\alpha)$, it follows that
\[
\delta\widetilde V_{\ext}(e,\alpha)
=
\left\langle
\grad_1\widetilde V_{\ext}(e,\alpha),
\eta
\right\rangle.
\]
Therefore,
\[
\begin{aligned}
0
=
\delta\widetilde{\mathcal A}_G[g]
&=
\int_a^b
\left\langle
\mathcal H(\xi),\delta\xi
\right\rangle
-
\left\langle
\grad_1\widetilde V_{\ext}(e,\alpha),
\eta
\right\rangle
dt\\
&=
\int_a^b
\left\langle
\mathcal H(\xi),\dot\eta
\right\rangle
+
\left\langle
\ad^\dagger_\xi\mathcal H(\xi)
-
\grad_1\widetilde V_{\ext}(e,\alpha),
\eta
\right\rangle
dt.
\end{aligned}
\]

Observe that $\grad\widetilde V$ is horizontal. Indeed, if $h\in G$ and
$X\in\Versp_h$, then
\[
\left\langle
\grad\widetilde V(h),X
\right\rangle_G
=
d\widetilde V_h(X)
=
dV_{\pi(h)}\bigl(T_h\pi(X)\bigr)
=
0.
\]
Moreover, if $Y_h\in\Horsp_h$ and
$Y:=T_h\pi(Y_h)$, then
\[
\left\langle
\grad\widetilde V(h),Y_h
\right\rangle_G
=
dV_{\pi(h)}(Y)
=
\left\langle
\grad V(\pi(h)),Y
\right\rangle_H.
\]
Hence $\grad\widetilde V$ is the horizontal lift of $\grad V$. In particular,
\[
T_hL_{h^{-1}}\bigl(\grad\widetilde V(h)\bigr)\in\mathfrak h.
\]
By the partial-symmetry identity,
\[
T_gL_{g^{-1}}\bigl(\grad\widetilde V(g)\bigr)
=
\grad_1\widetilde V_{\ext}(e,\alpha),
\]
and therefore $\grad_1\widetilde V_{\ext}(e,\alpha)\in\mathfrak h$.

Hence, by the definition of the distributional derivative,
\[
\frac{d}{dt}\mathcal H(\xi)
=
\ad^\dagger_\xi\mathcal H(\xi)
-
\grad_1\widetilde V_{\ext}(e,\alpha)
\]
in the weak sense.
Since $\mathcal H:\mathfrak g\to\mathfrak h$ is a fixed linear projection, this is precisely
Equation~\eqref{eq: EP eqs horizontal geodesics} in the weak sense.

Since $\xi\in L^2([a,b],\mathfrak g)$ and $\alpha\in H^1([a,b],M)$, the smoothness of the
potential implies that
$\grad_1\widetilde V_{\ext}(e,\alpha)$ is continuous. Moreover,
$\ad^\dagger_\xi\mathcal H(\xi)\in L^1([a,b],\mathfrak g)$. Hence
$\mathcal H(\xi)\in W^{1,1}([a,b],\mathfrak h)\subset C^0([a,b],\mathfrak h)$. Since $\mathcal H(\xi)$ is now continuous and hence bounded,
$\ad^\dagger_\xi\mathcal H(\xi)\in L^2([a,b],\mathfrak g)$, and therefore
$\mathcal H(\xi)\in H^1([a,b],\mathfrak h)$.

Finally, differentiating $\alpha=\Psi_{g^{-1}}(\alpha_0)$ in the weak sense gives $\dot\alpha=-\xi_M(\alpha)$, which proves Equation~\eqref{eq: parameter dynamics horizontal geodesics}.
\end{proof}

Observe that if $g$ is horizontal, then $\xi(t)\in\mathfrak h$ for almost every
$t\in[a,b]$, and Equation~\eqref{eq: EP eqs horizontal geodesics} is equivalent to
\[
\dot\xi
=
\ad^\dagger_{\xi}\xi
-
\grad_1\widetilde V_{\ext}(e,\alpha).
\]
These are the ordinary Euler-Poincar\'e equations for the corresponding mechanical
action on $G$. From Theorem
\ref{prop: variational principle equivalence homogeneous space geodesics},
$q=\pi\circ g$ is a critical point of
$\mathcal A_H$ and is therefore smooth. Since $g$ is horizontal, it
is the horizontal lift of $q$ with prescribed initial value, and is
therefore smooth as well. The difference, of course, is that non-horizontal solutions to the
ordinary Euler-Poincar\'e equations on $G$ need not project to critical points of
$\mathcal A_H$ in general. On the other hand, every critical point of
$\widetilde{\mathcal A}_G$ is obtained from a smooth horizontal lift by an
$H^1$ curve with values in $K$, as made precise in the following corollary.

\begin{corollary}\label{cor: horizontal_gauge}
Suppose that $\widetilde V$ admits a partial symmetry over some smooth manifold
$M$, and $g\in\Omega(G)$ is a critical point of
$\widetilde{\mathcal A}_G$.
Then, for every $h_a\in\pi^{-1}(\{\pi(g(a))\})$, there exists a unique horizontal curve
$h\in C^\infty([a,b],G)$ with $h(a)=h_a$ and a unique curve
$k\in H^1([a,b],K)$ such that $g(t)=h(t)k(t)$ for all $t\in[a,b]$. In such a case,
$\eta:=h^{-1}\dot h$ and $\beta:=\Psi_{h^{-1}}(\alpha_0)$ satisfy the Euler-Poincar\'e equations
\begin{align*}
\dot{\eta}
&=
\ad^\dagger_\eta\eta
-
\grad_1\widetilde V_{\ext}(e,\beta),\\
\dot{\beta}
&=
-\eta_M(\beta).
\end{align*}
Moreover, every critical point of $\widetilde{\mathcal A}_G$ can be recovered in
this way.
\end{corollary}

\begin{proof}
Let $g\in\Omega(G)$ be a critical point of
$\widetilde{\mathcal A}_G$, and let
$h_a\in\pi^{-1}(\{\pi(g(a))\})$. From Theorem
\ref{prop: variational principle equivalence homogeneous space geodesics},
the curve $q=\pi\circ g$ is a critical point of $\mathcal A_H$, and is therefore smooth. Let $h$ be the
unique horizontal lift of $q$ satisfying $h(a)=h_a$. Then $h$ is smooth and,
viewed on the corresponding fixed-endpoint path space, again by Theorem
\ref{prop: variational principle equivalence homogeneous space geodesics},
is a critical point of $\widetilde{\mathcal A}_G$.

Since $h$ is horizontal, Theorem~\ref{thm: EP_eqs_homo_G} implies that
$\eta:=h^{-1}\dot h$ and $\beta:=\Psi_{h^{-1}}(\alpha_0)$ satisfy
\begin{align*}
\dot{\eta}
&=
\ad^\dagger_\eta\eta
-
\grad_1\widetilde V_{\ext}(e,\beta),\\
\dot{\beta}
&=
-\eta_M(\beta).
\end{align*}

Setting $k:=h^{-1}g$, we have $k\in H^1([a,b],G)$. Since
$\pi\circ g=\pi\circ h$, it follows that
$k(t)=h(t)^{-1}g(t)\in K$ for every $t\in[a,b]$.
Hence $k\in H^1([a,b],K)$ and $g=hk$. The uniqueness of $k$ follows immediately from the uniqueness of the horizontal
lift $h$ with prescribed initial value.
\end{proof}

\subsection{Gauge invariance and optimal gauges}

The lifted functional $\widetilde{\mathcal A}_G$ possesses a natural gauge symmetry associated with the
isotropy subgroup $K$. The expression in~\eqref{eq: action function G geodesics}
extends naturally to a functional on $H^1([a,b],G)$, which we denote by the same
symbol $\widetilde{\mathcal A}_G$.

\begin{definition}\label{def: gauge_invariance}
The lifted action $\widetilde{\mathcal A}_G$ is said to be \textit{gauge invariant} under
$K$ if
\[
\widetilde{\mathcal A}_G[gk]
=
\widetilde{\mathcal A}_G[g]
\]
for every $g\in H^1([a,b],G)$ and $k\in H^1([a,b],K)$. When working on the
fixed-endpoint path space $\Omega(G)$, an endpoint-preserving gauge
transformation additionally satisfies $k(a)=k(b)=e$.
\end{definition}

The functional \eqref{eq: action function G geodesics} is gauge invariant in this sense. Indeed,
$\pi(g(t)k(t))=\pi(g(t))$ for all $t\in[a,b]$, and hence
\[
\widetilde{\mathcal A}_G[gk]
=
\mathcal A_H[\pi\circ(gk)]
=
\mathcal A_H[\pi\circ g]
=
\widetilde{\mathcal A}_G[g].
\]

Thus $\widetilde{\mathcal A}_G$ is completely insensitive to
vertical $K$-valued gauge transformations. This degeneracy motivates
the introduction of an additional variational criterion for selecting
a distinguished lift.

We denote the kinetic-energy functionals on $G$ and $H$ by
\[
E_G[g]:=\frac12\int_a^b\|\dot g(t)\|_G^2\,dt,
\qquad
E_H[q]:=\frac12\int_a^b\|\dot q(t)\|_H^2\,dt.
\]

To select a distinguished representative, we introduce an additional variational problem in the vertical directions. By Corollary~\ref{cor: horizontal_gauge}, every critical point of $\widetilde{\mathcal A}_G$ can be written uniquely, once the initial horizontal lift is fixed, as
$g=hk$, where $h\in C^\infty([a,b],G)$ is horizontal and $k\in H^1([a,b],K)$.

We equip $K$ with the Riemannian metric induced by the inclusion
$\iota:K\hookrightarrow G$, and denote the corresponding pointwise
norm and Riemannian distance by $\|\cdot\|_K$ and $d_K$,
respectively. Since this metric is left-invariant, each connected component of
$K$ is complete. We choose the gauge by minimizing its vertical kinetic energy,
\[
E_K[k]
=
\frac12\int_a^b\|\dot k(t)\|_K^2\,dt,
\]
among curves in $K$ with the endpoint values prescribed by the boundary conditions on $g$ and the chosen horizontal lift $h$.

Any minimizer is a smooth length-minimizing geodesic on $K$.
Writing $\sigma:=k^{-1}\dot k$, it satisfies
\[
\dot\sigma=\ad^{\dagger,K}_{\sigma}\sigma,
\qquad
\dot k=T_eL_k(\sigma),
\]
where $\ad^{\dagger,K}$ denotes the metric adjoint of the Lie bracket
on $\mathfrak s$ with respect to the induced metric on $K$.

This leads to the following result.

\begin{theorem}\label{thm: full_optimization_G}
Let $q\in\Omega(H)$ be a critical point of $\mathcal A_H$, and let
$h$ be the horizontal lift of $q$ satisfying $h(a)=g_a$. Suppose that $h(b)^{-1}g_b\in K_0$, where $K_0$ denotes the connected component of $K$ containing the
identity.

Then a curve $g\in\Omega(G)$ solves
\[
\min_{g\in\Omega(G)} E_G[g],
\qquad
\text{subject to }\pi\circ g=q,
\]
if and only if $g=hk$, where $k\in C^\infty([a,b],K)$ is a length-minimizing geodesic
satisfying $k(a)=e$ and $k(b)=h(b)^{-1}g_b$. In particular, $g$ is smooth and
\[
E_G[g]
=
E_H[q]
+
\frac{d_K\left(e,h(b)^{-1}g_b\right)^2}{2(b-a)}.
\]

Writing
$\sigma:=k^{-1}\dot k$, the minimizing gauge satisfies
\[
\dot\sigma=\ad^{\dagger,K}_{\sigma}\sigma,
\qquad
\sigma(t)\in\mathfrak s .
\]

Moreover, if $\widetilde V$ admits a partial symmetry over some
smooth manifold $M$, then, defining $\eta:=h^{-1}\dot h$ and $\beta:=\Psi_{h^{-1}}(\alpha_0)$, we have
\[
\dot\eta
=
\ad^\dagger_\eta\eta
-
\grad_1\widetilde V_{\ext}(e,\beta),
\qquad
\dot\beta=-\eta_M(\beta),
\]
with $\eta(t)\in\mathfrak h$, and the boundary conditions
\[
h(a)=g_a,\qquad
\pi(h(b))=q_b,\qquad
\beta(a)=\Psi_{g_a^{-1}}(\alpha_0).
\]
\end{theorem}

\begin{proof}
Since $q$ is a critical point of $\mathcal A_H$, it is smooth.
Let $h$ be its horizontal lift satisfying $h(a)=g_a$.

Let $g\in\Omega(G)$ be any curve satisfying $\pi\circ g=q$.
Since $\pi(g(t))=\pi(h(t))$ for every $t\in[a,b]$, the curve
$k(t):=h(t)^{-1}g(t)$ takes values in $K$. Moreover,
$k\in H^1([a,b],K)$, $k(a)=e$, and
$k(b)=h(b)^{-1}g_b$. Conversely, every
$k\in H^1([a,b],K)$ satisfying these endpoint conditions determines
an admissible lift $g=hk$ of $q$. Thus the admissible lifts of $q$
are in one-to-one correspondence with such curves $k$.

Let $\eta:=h^{-1}\dot h$, $\xi:=g^{-1}\dot g$, and
$\sigma:=k^{-1}\dot k$. Differentiating $g=hk$ and
left-trivializing gives
\[
\xi=\Ad_{k^{-1}}\eta+\sigma.
\]
Since $k(t)\in K$ and $\mathfrak h$ is $\Ad_K$-invariant,
$\Ad_{k^{-1}}\eta\in\mathfrak h$. Moreover,
$\sigma\in\mathfrak s$, and
$\mathfrak g=\mathfrak s\oplus\mathfrak h$ is an orthogonal
decomposition. Hence
\[
\|\xi\|_G^2
=
\|\Ad_{k^{-1}}\eta\|_G^2+\|\sigma\|_G^2
=
\|\eta\|_G^2+\|\sigma\|_G^2,
\]
where the second equality follows from the $\Ad_K$-invariance of the
metric on $\mathfrak h$. By left-invariance of the induced metric on
$K$, $\|\sigma\|_G=\|\dot k\|_K$. Therefore,
\[
E_G[g]
=
E_G[h]+E_K[k]
=
E_H[q]+E_K[k],
\]
where the last equality follows because $h$ is horizontal and
$\pi:G\to H$ is a Riemannian submersion.

Thus minimizing $E_G[g]$ among all lifts of the fixed curve $q$
is equivalent to minimizing $E_K[k]$ among curves in $K$ satisfying
$k(a)=e$ and $k(b)=h(b)^{-1}g_b$. The assumption
$h(b)^{-1}g_b\in K_0$ ensures that the endpoints lie in the same
connected component. Since $K_0$ is complete, the Hopf--Rinow
theorem guarantees the existence of a length-minimizing geodesic
joining them. The minimizers of $E_K$ are precisely the
constant-speed length-minimizing geodesics with these endpoints.
Hence $k$ is smooth and, with $\sigma=k^{-1}\dot k$,
\[
\dot\sigma=\ad_{\sigma}^{\dagger,K}\sigma.
\]
Since $h$ and $k$ are smooth, so is $g=hk$.

Finally, a length-minimizing geodesic has constant speed
\[
\|\dot k(t)\|_K
=
\frac{d_K\bigl(e,h(b)^{-1}g_b\bigr)}{b-a},
\]
and therefore
\[
E_K[k]
=
\frac{d_K\bigl(e,h(b)^{-1}g_b\bigr)^2}{2(b-a)}.
\]
Substitution into $E_G[g]=E_H[q]+E_K[k]$ gives the stated
energy formula.

If, in addition, $\widetilde V$ admits a partial symmetry over some
smooth manifold $M$, then, on the fixed-endpoint path space with
endpoints $g_a$ and $h(b)$, Theorem
\ref{prop: variational principle equivalence homogeneous space geodesics}
implies that $h$ is a critical point of
$\widetilde{\mathcal A}_G$. Since $h$ is horizontal,
Theorem~\ref{thm: EP_eqs_homo_G} gives, with
$\beta:=\Psi_{h^{-1}}(\alpha_0)$,
\[
\dot\eta
=
\ad^\dagger_\eta\eta
-
\grad_1\widetilde V_{\ext}(e,\beta),
\qquad
\dot\beta=-\eta_M(\beta).
\]
The stated boundary conditions follow from the definitions of
$h$ and $\beta$.
\end{proof}

\begin{remark}
Suppose that $q$ is a closed curve, so that $q(a)=q(b)$, and let $h$ be its horizontal lift with $h(a)=g_a$. Since $h(a)$ and $h(b)$ belong to the same fiber, there exists a unique $\kappa_q\in K$ such that $h(b)=g_a\kappa_q$. The element $\kappa_q$ is the holonomy element determined by the horizontal lift of $q$ based at $g_a$.

If we additionally require the lifted curve to be closed, $g(a)=g(b)=g_a$, then the vertical gauge $k$ in the decomposition $g=hk$ must satisfy
$k(a)=e$ and $k(b)=\kappa_q^{-1}$. Provided that $\kappa_q\in K_0$, Theorem~\ref{thm: full_optimization_G} applies. Moreover,
by symmetry and left-invariance of $d_K$,
$d_K(e,\kappa_q^{-1})=d_K(e,\kappa_q)$. Therefore,
\[
\min_{\substack{g(a)=g(b)=g_a\\ \pi\circ g=q}}
E_G[g]
=
E_H[q]
+
\frac{d_K(e,\kappa_q)^2}{2(b-a)}.
\]
Thus the additional energy is precisely the minimum vertical kinetic energy required to compensate for the holonomy of the horizontal lift and close the curve in $G$. In particular, this contribution vanishes exactly when the horizontal lift closes, i.e., when $\kappa_q=e$.
\end{remark}

\subsection{Example: Pure quantum states and unitary gauges}

Let $\mathbb C^{n+1}$ be the Hilbert space of an $(n+1)$-level quantum
system. A pure quantum state is represented by a ray $[\psi]$, since
normalized vectors differing only by a complex phase describe the same
physical state. Equivalently, using bra-ket notation, $[\psi]$ may be represented by the
rank-one orthogonal projector
$\rho=|\psi\rangle\langle\psi|$, characterized by
$\rho^\dagger=\rho$, $\rho^2=\rho$, and $\tr(\rho)=1$ (see \cite{BengtssonZyczkowski} for more details). Hence
\[
\CP^n\cong
\mathcal P_1:=
\left\{
\rho\in\operatorname{End}(\mathbb C^{n+1}) \, \vert \,
\rho^\dagger=\rho,\ \rho^2=\rho,\ \tr(\rho)=1
\right\}.
\]

The group $G=SU(n+1)$ acts transitively on $\mathcal P_1$ by
conjugation, i.e. $U\cdot\rho=U\rho U^\dagger$. Fix a reference state
$\rho_0=|\psi_0\rangle\langle\psi_0|$. Its isotropy subgroup is
$$K=\{U\in SU(n+1):U\rho_0U^\dagger=\rho_0\}
\cong S(U(1)\times U(n)),$$
and the quotient map is
$$\pi:SU(n+1)\to\CP^n, \qquad \pi(U) = U\rho_0U^\dagger,$$
so that $\CP^n\cong SU(n+1)/K$. A curve $U(t)\in SU(n+1)$ is
therefore a unitary lift of the physical state trajectory
$\rho(t)=U(t)\rho_0U(t)^\dagger$. This lift is not unique, since $U$ and
$Uk$ determine the same trajectory $\rho(t)$ for any $k(t)\in K$.

We equip $\su(n+1)$ with the $\Ad$-invariant inner product
$\langle A,B\rangle_G=-\tr(AB)$. The vertical Lie algebra is
$\mathfrak s=\{A\in\su(n+1):[A,\rho_0]=0\}$, and
$\mathfrak h=\mathfrak s^\perp$.
Thus elements of $\mathfrak s$ generate unitary motions that leave the
physical state unchanged, whereas elements of $\mathfrak h$ generate motion in projective Hilbert space. We denote by
$\mathcal H:\su(n+1)\to\mathfrak h$ the corresponding orthogonal
projection. With the convention
$g_{\mathrm{FS}}(\delta\rho_1,\delta\rho_2)
=\frac12\tr(\delta\rho_1\delta\rho_2)$, the quotient metric induced by
$\langle\cdot,\cdot\rangle_G$ satisfies
\[
\langle\delta\rho_1,\delta\rho_2\rangle_H
=
\tr(\delta\rho_1\delta\rho_2)
=
2g_{\mathrm{FS}}(\delta\rho_1,\delta\rho_2).
\]
Hence the quotient geometry considered in Section~3 is precisely the
Fubini--Study geometry, up to this fixed normalization \cite{BengtssonZyczkowski, helgason2024geometric}.

Let $\rho_d=|\psi_d\rangle\langle\psi_d|$ be a desired pure state and
consider the potential
\[
V(\rho)
=
\lambda\bigl(1-\tr(\rho\rho_d)\bigr),
\qquad \lambda>0.
\]
For pure states, the quantity
$F(\rho,\rho_d):=\tr(\rho\rho_d)
=|\langle\psi,\psi_d\rangle|^2\in[0,1]$
is the fidelity between the two states, i.e., their squared quantum-state
overlap. It equals one when they represent the same physical state and
zero when they are orthogonal. Thus
$V(\rho)=\lambda(1-F(\rho,\rho_d))$ is a scaled infidelity that
penalizes deviation from the desired state. Its lift to $SU(n+1)$ is
$$\widetilde V(U)
=\lambda(1-\tr(U\rho_0U^\dagger\rho_d)).$$
Take $M=\mathcal P_1$ with conjugation action
$\Psi_W(\Gamma)=W\Gamma W^\dagger$ and define
$\widetilde V_{\mathrm{ext}}(U,\Gamma)
=
\lambda\left(
1-\tr(U\rho_0U^\dagger\Gamma)
\right)$. Then
$\widetilde V_{\mathrm{ext}}(U,\rho_d)=\widetilde V(U)$ and
$\widetilde V_{\mathrm{ext}}(WU,W\Gamma W^\dagger)
=\widetilde V_{\mathrm{ext}}(U,\Gamma)$, so that
$\widetilde V$ admits a partial symmetry.

Let $\xi=U^{-1}\dot U$ and define
$\alpha=U^\dagger\rho_dU$. Thus $\alpha$ is the desired state
expressed in the moving unitary frame defined by $U(t)$. The
infinitesimal conjugation action gives
$\dot\alpha=-[\xi,\alpha]$. A direct calculation yields
$\grad_1\widetilde V_{\mathrm{ext}}(e,\alpha)
=
\lambda[\rho_0,\alpha]\in\mathfrak h$. Since the trace metric is $\Ad$-invariant, it is easily see that
$\ad^\dagger_\xi\eta=[\eta,\xi]$, and hence Theorem~\ref{thm: EP_eqs_homo_G}  yields
\[
\frac{d}{dt}\mathcal H(\xi)
=
[\mathcal H(\xi),\xi]
-\lambda[\rho_0,\alpha],
\qquad
\dot\alpha=-[\xi,\alpha].
\]
The horizontal component $\mathcal H(\xi)$ changes the physical state,
whereas the vertical component
$\xi_{\mathfrak s}:=\xi-\mathcal H(\xi)$ changes only its unitary
representative. In particular, if $h$ is the horizontal representative
and
$\eta=h^{-1}\dot h\in\mathfrak h$,
$\beta=h^\dagger\rho_dh$, then
$\dot\eta=-\lambda[\rho_0,\beta]$ and
$\dot\beta=-[\eta,\beta]$.

This decomposition has a direct quantum-control interpretation.
A finite-dimensional controlled quantum system is typically generated
by a Hamiltonian of the form
\[
\mathsf H(t)
=
\mathsf H_0+\sum_{j=1}^m u_j(t)\mathsf H_j,
\qquad
i\dot U(t)=\mathsf H(t)U(t),
\]
where, after removing an irrelevant scalar multiple of the identity,
the Hamiltonian may be taken traceless when working on $SU(n+1)$.
The available control Hamiltonians determine the admissible unitary
motions and the associated controllability properties; see, for
instance,~\cite{Altafini2002}. For every admissible propagator,
$\xi=-iU^\dagger\mathsf H U$. Hence the splitting
$\xi=\xi_{\mathfrak h}+\xi_{\mathfrak s}$ separates the
body-represented Hamiltonian into a state-changing component and an
isotropy component. Indeed, $\dot\rho=U[\xi,\rho_0]U^\dagger$, while
$[\xi_{\mathfrak s},\rho_0]=0$.

The optimal-gauge problem therefore acquires a concrete control
interpretation. Fix a critical physical state trajectory $\rho(t)$ obtained from the
variational problem above and prescribe unitary representatives
$U_a,U_b\in SU(n+1)$ at the endpoints. Let
$h$ be the horizontal lift satisfying $h(a)=U_a$ and define the
endpoint mismatch $\kappa:=h(b)^{-1}U_b\in K$.
Since $K\cong S(U(1)\times U(n))$ is connected and the trace metric
restricts to a bi-invariant metric on $K$, choose
$S_*\in\mathfrak s$ of minimum norm among the logarithms of $\kappa$,
that is,
$\Exp(S_*)=\kappa$ and
$\|S_*\|=\min\{\|S\|:\Exp(S)=\kappa\}$.
Theorem~\ref{thm: full_optimization_G} then gives the optimal gauge and lift explicitly:
\[
k_*(t)
=
\Exp\left(
\frac{t-a}{b-a}S_*
\right),
\qquad
U_*(t)=h(t)k_*(t).
\]
The optimal gauge has constant body generator
$\sigma_{*}=S_{*} /(b-a)$ and satisfies
\[
E_G[U_{*}]
=
E_H[\rho]
+
\frac{d_K(e,\kappa)^2}{2(b-a)}.
\]
Thus $U_*$ is a closed-form solution of the optimal lift problem:
among all unitary evolutions producing the prescribed pure-state
trajectory and satisfying the prescribed unitary endpoints, it
minimizes the kinetic energy induced by the trace metric.

The corresponding Hamiltonian also separates into horizontal and
gauge components. Since $U_*=hk_*$ and
$\eta=h^{-1}\dot h$, one obtains
\[
\mathsf H_*(t)
=
i\,h(t)\eta(t)h(t)^\dagger
+
\mathsf H_{\mathrm{gauge}}^*(t),
\qquad
\mathsf H_{\mathrm{gauge}}^*(t)
=
i\,h(t)\sigma_*h(t)^\dagger.
\]
Moreover, $[\mathsf H_{\mathrm{gauge}}^*(t),\rho(t)]=0$, because $\sigma_*\in\mathfrak s$. Hence the optimal gauge Hamiltonian
does not alter the prescribed physical state trajectory; it is exactly
the minimum additional Hamiltonian motion required to satisfy the
terminal unitary representative. For traceless Hamiltonians,
\[
E_G[U]
=
\frac12\int_a^b\tr\bigl(\mathsf H(t)^2\bigr)\,dt,
\]
so the energy splitting separates the generator effort required to
move the physical state from the minimum additional effort spent
entirely in isotropy directions. If the terminal unitary
representative is free, the latter contribution vanishes and the
horizontal lift itself is the minimum-energy representative.

This interpretation is exact whenever the optimal Hamiltonian is
realizable by the available controls and the control cost agrees with
the trace metric. For a restricted system
$\mathsf H=\mathsf H_0+\sum_j u_j\mathsf H_j$, endpoint
controllability alone does not imply that
$\mathsf H_{\mathrm{gauge}}^*(t)$ is instantaneously admissible. In that case, for the same trace-metric cost and prescribed state
trajectory, the construction gives the intrinsic geometric optimum,
and therefore a natural lower bound and benchmark for the corresponding
constrained control problem.

For a two-level system,
$\CP^1\cong SU(2)/U(1)\cong S^2$. Writing a pure state as
$\rho=\frac12(I+r\cdot\boldsymbol{\sigma})$, where
$r\in S^2$ and $\boldsymbol{\sigma}$ denotes the vector of Pauli
matrices, the fidelity potential becomes
$V(r)
=
\frac{\lambda}{2}(1-r\cdot r_d)$. Thus the fidelity potential penalizes misalignment of the Bloch vector
with the desired state. Since the isotropy subgroup is $U(1)$, the optimal vertical
generator is constant and produces the least costly rotation about the
instantaneous Bloch direction required to realize the prescribed
unitary endpoint, without changing the Bloch vector itself. For a
closed Bloch trajectory, the same correction compensates the holonomy
of the horizontal lift and gives the minimum additional unitary energy
required to close the lift.

\section{Application: Optimal Camera Orientation on a Satellite with an Undesirable Pointing Region}
\label{sec:satellite_application}

Satellites can be modeled as rigid bodies on $G=\SO(3)$. We consider a camera fixed to a satellite with pointing direction $q=Re_3\in S^2$ in the inertial frame, where $R\in\SO(3)$ is the satellite orientation and $e_3=(0,0,1)^T\in\R^3$ is the camera direction in the body frame. Let $K\subset\SO(3)$ be the Lie subgroup given by the elements of the form 
\[ k(\theta) = \begin{bmatrix} \cos(\theta) & \sin(\theta) & 0 \\ -\sin(\theta) & \cos(\theta) & 0 \\ 0 & 0 & 1 \end{bmatrix}, \qquad \theta\in\R/2\pi\Z . \] 
It is clear that $K\cong\SO(2)$, and that $K=\operatorname{Stab}(e_3)$ with respect to the left-$\SO(3)$ action on $S^2$ defined by 
\[ \Phi:\SO(3)\times S^2\to S^2, \qquad (R,q)\mapsto Rq . \] 
Moreover, $\Phi$ is transitive, so that the quotient space $H=\SO(3)/\SO(2)\cong S^2$ is a homogeneous space of $\SO(3)$, and the quotient map $\pi:\SO(3)\to S^2$ is given by $\pi(R)=Re_3$. Let $\hat{\cdot}:\R^3\to\so(3)$ denote the \emph{hat isomorphism}, and let $\cdot^\vee$ denote its inverse, where 
\[ \widehat{ \begin{bmatrix} x \\ y \\ z \end{bmatrix}} = \begin{bmatrix} 0 & -z & y \\ z & 0 & -x \\ -y & x & 0 \end{bmatrix}. \] 
We equip $\SO(3)$ with the left-invariant metric determined at the identity by 
\[ \left<\hat\Omega_1,\hat\Omega_2\right>_{\SO(3)} = \Omega_1^T\mathbb J\Omega_2, \] 
for all $\Omega_1,\Omega_2\in\R^3$, where $\mathbb J\in\operatorname{Sym}_3^+(\R)$ is the inertia tensor of the satellite. We take the origin of the body frame to be the satellite's center of mass, and the coordinate axes of the body to be the principal axes of inertia, so that 
\[ \mathbb J=\operatorname{diag}(J_1,J_2,J_3). \] 
Consider the basis $\{\hat\Omega_1,\hat\Omega_2,\hat\Omega_3\}$ of $\mathfrak g=\so(3)$ given by 
\begin{align}\label{eq: orthonormal basis S^2} &\hat\Omega_1 = \begin{bmatrix} 0 & 0 & 0 \\ 0 & 0 & -1 \\ 0 & 1 & 0 \end{bmatrix}, && \hat\Omega_2 = \begin{bmatrix} 0 & 0 & 1 \\ 0 & 0 & 0 \\ -1 & 0 & 0 \end{bmatrix}, && \hat\Omega_3 = \begin{bmatrix} 0 & -1 & 0 \\ 1 & 0 & 0 \\ 0 & 0 & 0 \end{bmatrix}. 
\end{align} 
Thus $\hat\Omega_1=\hat e_1$, $\hat\Omega_2=\hat e_2$, and $\hat\Omega_3=\hat e_3$. Moreover, \[ T_I\pi(\hat\Omega)=\hat\Omega e_3=\Omega\times e_3, \] so that \[ T_I\pi(\hat\Omega_1)=-e_2, \qquad T_I\pi(\hat\Omega_2)=e_1, \qquad T_I\pi(\hat\Omega_3)=0. \] Hence 
\[ \mathfrak s=\ker(T_I\pi)=\operatorname{span}\{\hat\Omega_3\}. \]
It is easily verified that $\langle\hat\Omega_i,\hat\Omega_j\rangle_{\SO(3)}=0$ whenever $i\neq j$, so that
\[ \mathfrak h=\mathfrak s^\perp = \operatorname{span}\{\hat\Omega_1,\hat\Omega_2\}. \]
For the quotient construction of Section~\ref{sec6}, the inertia inner product on $\so(3)$ must be $\Ad_K$-invariant. For 
\[ k(\theta) = \begin{bmatrix} \cos(\theta) & \sin(\theta) & 0 \\ -\sin(\theta) & \cos(\theta) & 0 \\ 0 & 0 & 1 \end{bmatrix}, \]
we find by direct calculation that 
\begin{align*} 
\Ad_{k(\theta)}\hat\Omega_1 &= \cos(\theta)\hat\Omega_1-\sin(\theta)\hat\Omega_2,\\ \Ad_{k(\theta)}\hat\Omega_2 &= \sin(\theta)\hat\Omega_1+\cos(\theta)\hat\Omega_2,\\
\Ad_{k(\theta)}\hat\Omega_3 &= \hat\Omega_3 . 
\end{align*} 
Hence, 
\begin{align*} \left< \Ad_{k(\theta)}\hat\Omega_1, \Ad_{k(\theta)}\hat\Omega_1 \right>_{\SO(3)} &= \cos^2(\theta)J_1+\sin^2(\theta)J_2,\\
\left< \Ad_{k(\theta)}\hat\Omega_2, \Ad_{k(\theta)}\hat\Omega_2 \right>_{\SO(3)} &= \sin^2(\theta)J_1+\cos^2(\theta)J_2,\\
\left< \Ad_{k(\theta)}\hat\Omega_1, \Ad_{k(\theta)}\hat\Omega_2 \right>_{\SO(3)} &= \sin(\theta)\cos(\theta)(J_1-J_2). 
\end{align*} 
Therefore, the inertia inner product is $\Ad_K$-invariant if and only if $J_1=J_2=:J$. We assume this condition throughout. Thus the satellite is axisymmetric about $e_3$, and the camera axis is aligned with its symmetry axis. In such a case, the induced inner product on $T_{e_3}S^2$ is given by \[ \left<X,Y\right>_{T_{e_3}S^2} := \left< \left(T_I\pi\big|_{\mathfrak h}\right)^{-1}X, \left(T_I\pi\big|_{\mathfrak h}\right)^{-1}Y \right>_{\SO(3)} . \] Since $\{e_1,e_2\}$ forms a basis for $T_{e_3}S^2$, we set $X=X^1e_1+X^2e_2$ and $Y=Y^1e_1+Y^2e_2$. Since \[ \left(T_I\pi\big|_{\mathfrak h}\right)^{-1}e_1=\hat\Omega_2, \qquad \left(T_I\pi\big|_{\mathfrak h}\right)^{-1}e_2=-\hat\Omega_1, \] we calculate \begin{align*} \left<X,Y\right>_{T_{e_3}S^2} &= \left< X^1\hat\Omega_2-X^2\hat\Omega_1, Y^1\hat\Omega_2-Y^2\hat\Omega_1 \right>_{\SO(3)}\\ &= JX^1Y^1+JX^2Y^2\\ &= JX\cdot Y. \end{align*} We extend this inner product to an $\SO(3)$-invariant Riemannian metric on $S^2$ by the left action. That is, \begin{equation}\label{eq: metric S^2} \left<X,Y\right>_{S^2} := \left<R^TX,R^TY\right>_{T_{e_3}S^2} = J(R^TX)\cdot(R^TY) = JX\cdot Y \end{equation} for all $X,Y\in T_qS^2$ and for any $R\in\SO(3)$ such that $\pi(R)=q$. Hence $\langle\cdot,\cdot\rangle_{S^2}$ is simply a scalar multiple of the round metric, and $\pi:\SO(3)\to S^2$ is a Riemannian submersion. Moreover, the commutator relations \[ [\hat\Omega_1,\hat\Omega_2]=\hat\Omega_3, \qquad [\hat\Omega_1,\hat\Omega_3]=-\hat\Omega_2, \qquad [\hat\Omega_2,\hat\Omega_3]=\hat\Omega_1 \] show that $\mathfrak g=\mathfrak s\oplus\mathfrak h$ satisfies \emph{Cartan's decomposition}:
\begin{equation}\label{cartan}
    [\mathfrak{s}, \mathfrak{s}] \subset \mathfrak{s}, \quad [\mathfrak{s}, \mathfrak{h}] \subset \mathfrak{h}, \quad [\mathfrak{h}, \mathfrak{h}] \subset \mathfrak{s}.
\end{equation}
Thus $S^2\cong\SO(3)/\SO(2)$ equipped with the above Riemannian metric is a Riemannian symmetric space (see \cite{helgason2024geometric} for more details). 

\subsection{Projected camera-pointing problem} Suppose that the satellite motion is governed by the fully actuated kinematic control law \[ \dot R=R\hat\Omega, \] where $\Omega\in\R^3$ is the body angular velocity. We consider the time interval $[0,1]$ and prescribe initial and final orientations $R_0,R_1\in\SO(3)$. Therefore, the camera must move from 
\[ q_0:=\pi(R_0)=R_0e_3 \] to \[ q_1:=\pi(R_1)=R_1e_3. \] 
We solve the task in two stages. First, we determine an optimal camera-pointing trajectory joining $q_0$ to $q_1$. Then, keeping this projected trajectory fixed, we use Theorem~\ref{thm: full_optimization_G} to select a minimum-energy lift satisfying the full attitude endpoint conditions. Let $\bar q\in S^2$ denote an undesirable pointing direction, such as the direction of the Sun, and define the geodesic cap 
\[ \mathcal C_r(\bar q) := \{q\in S^2:d_{S^2}(q,\bar q)<r\}, \] 
where $0<r<\pi\sqrt J$. We wish to construct $q(t)$ so that it connects prescribed boundary conditions, avoids the cap $\mathcal{C}_r(\bar q)$, and approximately minimizes energy. To do this, we use the strategy of artificial potential functions (see \cite{GoodmanThesis}). 

By~\eqref{eq: metric S^2}, 
\[ d_{S^2}(q,\bar q) = \sqrt J\,\arccos(q\cdot\bar q). \] 
We set $\rho_r:=1-\cos(r/\sqrt J),$ so that $d_{S^2}(q,\bar q)=r$ if and only if $1-q\cdot\bar q=\rho_r.$ Consider the smooth artificial potential \begin{equation}\label{eq:satellite_smooth_potential} 
V_{\rm obs}(q) := \frac{\tau} {1+(1-q\cdot\bar q)/\rho_r}, \qquad \tau>0,
\end{equation} 
which satisfies $V_{\rm obs}(\bar q)=\tau$, $V_{\rm obs}(q)=\tau/2$ on the boundary of $\mathcal C_r(\bar q)$, and decreases monotonically with geodesic distance from $\bar q$. This potential provides a smooth penalty rather than a hard state constraint; strict exclusion of $\mathcal C_r(\bar q)$ could instead be imposed separately. The control affects the camera direction via \[ \dot q = T_R\pi(\dot R) = R(\Omega\times e_3). \] Setting $u:=\dot q\in T_qS^2$, we consider the optimal control problem 
\begin{equation}\label{eq: optimal_control_S2} \begin{aligned} \min_{q,u}\quad& \int_0^1 \left( \frac12\|u\|_{S^2}^2+V_{\rm obs}(q) \right)dt,\\ \text{subject to}\quad& \dot q=u,\qquad q(0)=q_0,\qquad q(1)=q_1. \end{aligned} \end{equation}
Eliminating $u$, this is equivalent to minimizing \begin{equation}\label{eq:satellite_projected_action} \mathcal A_{S^2}[q] := \int_0^1 \left( \frac12\|\dot q\|_{S^2}^2+V_{\rm obs}(q) \right)dt \end{equation} 
over $q\in\Omega^{0,1}_{q_0,q_1}(S^2)$. To place this problem in the convention of Section~\ref{sec6}, set \[ W:=-V_{\rm obs}, \qquad \widetilde W:=W\circ\pi. \] The corresponding lifted action is \begin{equation}\label{eq:satellite_lifted_action} \widetilde{\mathcal A}_{\SO(3)}[R] = \int_0^1 \left( \frac12\|\mathcal H_R(\dot R)\|_{\SO(3)}^2 - \widetilde W(R) \right)dt. \end{equation} 
By Theorem~\ref{prop: variational principle equivalence homogeneous space geodesics}, critical points of~\eqref{eq:satellite_projected_action} are precisely the projections of critical points of~\eqref{eq:satellite_lifted_action}. For $\dot R=R\hat\Omega$, define \[ \Omega_{\mathfrak h} := \left( T_RL_{R^{-1}}(\mathcal H_R(\dot R)) \right)^\vee . \] Then \[ \Omega_{\mathfrak h} = e_3\times(\Omega\times e_3) = \Omega-(\Omega\cdot e_3)e_3, \] and \begin{equation}\label{eq: horizontal projection S^2 norm} \|\mathcal H_R(\dot R)\|_{\SO(3)}^2 = J(\Omega_1^2+\Omega_2^2) = J\|\Omega\times e_3\|^2 = \|\dot q\|_{S^2}^2. \end{equation} 

\subsection{Symmetry breaking and the horizontal representative} 
Let $q$ be a minimizer of~\eqref{eq: optimal_control_S2}, and hence a critical point of $\mathcal A_{S^2}$. Let $R_h$ be its horizontal lift satisfying $R_h(0)=R_0.$ Thus 
\[ R_h(t)e_3=q(t). \] 
To apply the reduced equations of Section~\ref{sec6}, take the parameter manifold to be $M=S^2$ with the standard left action $\Psi_P(a)=Pa$. 
Define 
\begin{equation}\label{eq:satellite_extended_potential} \widetilde W_{\ext}(P,a) := - \frac{\tau} {1+\bigl(1-(Pe_3)\cdot a\bigr)/\rho_r}. 
\end{equation} 
Then $\widetilde W_{\ext}(P,\bar q)=\widetilde W(P),$ and $\widetilde W_{\ext}(SP,Sa)=\widetilde W_{\ext}(P,a)$ for all $S,P\in\SO(3)$ and $a\in S^2$, since $S$ preserves the Euclidean inner product. Hence $\widetilde W$ admits a partial symmetry. Introduce the horizontal body velocity 
\[ \hat\Pi:=R_h^{-1}\dot R_h\in\mathfrak h \] 
and the advected parameter \[ \alpha:=R_h^{-1}\bar q\in S^2. \] Thus $\alpha$ is the undesirable pointing direction expressed in the moving body frame. Its evolution is \[ \dot\alpha=-\Pi\times\alpha. \] For the inertia metric, one has \[ \left(\ad^\dagger_{\hat\Omega}\hat\sigma\right)^\vee = \mathbb J^{-1}\bigl((\mathbb J\sigma)\times\Omega\bigr). \] Since $\Pi\cdot e_3=0$ and $J_1=J_2=J$, we have $\mathbb J\Pi=J\Pi$, and therefore $\ad^\dagger_{\hat\Pi}\hat\Pi=0$. Define \[ \chi(\alpha):=1+\frac{1-e_3\cdot\alpha}{\rho_r}. \] A direct differentiation of~\eqref{eq:satellite_extended_potential} gives \begin{equation}\label{eq:satellite_reduced_gradient} \grad_1\widetilde W_{\ext}(I,\alpha) = - \frac{\tau} {J\rho_r\chi(\alpha)^2} \widehat{e_3\times\alpha}. \end{equation} In particular, the reduced gradient belongs to $\mathfrak h$. The horizontal Euler--Poincar\'e equations and reconstruction equation therefore reduce to \begin{align} \dot\Pi &= \frac{\tau} {J\rho_r\chi(\alpha)^2} (e_3\times\alpha), \label{eq:satellite_reduced_Pi}\\ \dot\alpha &= -\Pi\times\alpha, \label{eq:satellite_reduced_A}\\ \dot R_h &= R_h\hat\Pi. \label{eq:satellite_horizontal_reconstruction} \end{align} The boundary conditions are \[ R_h(0)=R_0, \qquad R_h(1)e_3=q_1, \qquad \alpha(0)=R_0^{-1}\bar q, \] together with \[ \Pi(t)\cdot e_3=0, \qquad t\in[0,1]. \] Thus the projected problem determines the smooth horizontal representative $R_h$ and its horizontal body velocity $\Pi$. 

\subsection{Optimal vertical gauge} 
The projected trajectory fixes the camera direction, but not the complete satellite orientation. Indeed, every lift $R=R_h\kappa$, with $\kappa:[0,1]\to K\cong\SO(2)$, satisfies $\pi(R_h\kappa)=\pi(R_h)=q$. This is precisely the gauge freedom of Definition~\ref{def: gauge_invariance}. We now impose the prescribed terminal attitude $R_1$. Since $$R_h(1)e_3=q_1=R_1e_3,$$ 
the endpoint mismatch $\kappa_b:=R_h(1)^{-1}R_1$ belongs to $K$. Since $K\cong\SO(2)$ is connected, the hypothesis $\kappa_b\in K_0$ of Theorem~\ref{thm: full_optimization_G} is automatically satisfied. We therefore consider \begin{equation}\label{eq:satellite_optimal_lift} \min_{ R\in\Omega^{0,1}_{R_0,R_1}(\SO(3)) } E_{\SO(3)}[R], \qquad \text{subject to }\pi\circ R=q. \end{equation} By Theorem~\ref{thm: full_optimization_G}, the solution is $R^\ast=R_h\kappa^\ast$, where $\kappa^\ast$ is a length-minimizing geodesic in $K$ joining $I$ to $\kappa_b$. Choose a principal angle $\bar\theta\in[-\pi,\pi]$ such that $\kappa_b=\Exp(-\bar\theta\,\hat e_3)$. Then \begin{equation}\label{eq:satellite_optimal_gauge} \kappa^\ast(t) = \Exp(-t\bar\theta\,\hat e_3), \qquad t\in[0,1], \end{equation} and its constant body velocity is $\hat\Sigma^\ast:=(\kappa^\ast)^{-1}\dot\kappa^\ast =-\bar\theta\,\hat e_3$. If $|\bar\theta|=\pi$, there are two minimizing choices of principal angle, both giving the same minimum energy. The body angular velocity of the complete lift is \begin{equation}\label{eq:satellite_full_velocity} \hat\Omega^\ast = (R^\ast)^{-1}\dot R^\ast = \Ad_{(\kappa^\ast)^{-1}}\hat\Pi + \hat\Sigma^\ast . \end{equation} Thus $\hat\Pi$ determines the camera motion on $S^2$, whereas $\hat\Sigma^\ast$ is the minimum-energy rotation about the camera axis required to satisfy the prescribed terminal attitude. Since $\|\hat e_3\|_K^2=J_3$, we have $d_K(I,\kappa_b)=\sqrt{J_3}\,|\bar\theta|$, where \[ |\bar\theta| = \arccos \left( \frac{\operatorname{tr}(\kappa_b)-1}{2} \right). \] Hence $E_K[\kappa^\ast]=\frac{J_3}{2}\bar\theta^2$, and the energy splitting of Theorem~\ref{thm: full_optimization_G} becomes \begin{equation}\label{eq:satellite_energy_splitting} E_{\SO(3)}[R^\ast] = E_{S^2}[q] + \frac{J_3}{2}\bar\theta^2. \end{equation} The second term is therefore exactly the minimum rotational energy required to satisfy the full attitude endpoint condition without altering the prescribed camera trajectory.

\subsection{Numerical illustration} We solve the boundary-value problem~\eqref{eq:satellite_reduced_Pi}--\eqref{eq:satellite_horizontal_reconstruction} using a Lie group variational integrator (LGVI) on $\SO(3)$~\cite{LeeLeokMcClamroch2007}. In this subsection, $R_k$ denotes the discrete horizontal representative and $\Omega_k$ denotes its horizontal body angular velocity. For a step size $h=1/N$, the discrete update is \begin{align} h\widehat{\mathbb J\Omega_k} &= F_k\mathbb J_d-\mathbb J_dF_k^T, \nonumber\\ R_{k+1} &= R_kF_k, \nonumber\\ \mathbb J\Omega_{k+1} &= F_k^T\mathbb J\Omega_k+hM_{k+1}, \label{eq:satellite_LGVI} \end{align} where $\mathbb J_d=\frac12\tr(\mathbb J)I-\mathbb J$, $F_k\in\SO(3)$, and $M_{k+1}$ is the body moment induced by the symmetry-breaking potential. The endpoint condition $R_Ne_3=q_1$ is imposed by shooting. We take $J=1$, $J_3=1.6$, $N=240$, $\tau=3$, $r=30^\circ$, $\bar q=e_1$, and the time interval $[0,1]$. The numerical illustration has two separate purposes. First, we compare the projected trajectory obtained with the undesirable-pointing penalty to the unpenalized projected trajectory. This shows the effect of the symmetry-breaking potential on the camera path in $S^2$. Second, after the penalized projected trajectory has been fixed, we compare the minimum-energy vertical gauge from Theorem~\ref{thm: full_optimization_G} with a nonoptimal gauge having the same endpoint values. This second comparison illustrates the optimal-gauge theorem and does not change the projected camera trajectory\footnote{The code used to generate the simulations, as well as a supplementary video illustrating the projected camera-pointing trajectory, the evolution of the horizontal body frame, and the minimum-energy gauge reconstruction can be found at \url{https://doi.org/10.5281/zenodo.21996815}}. 

\begin{figure}[t] \centering \includegraphics[width=0.72\linewidth] {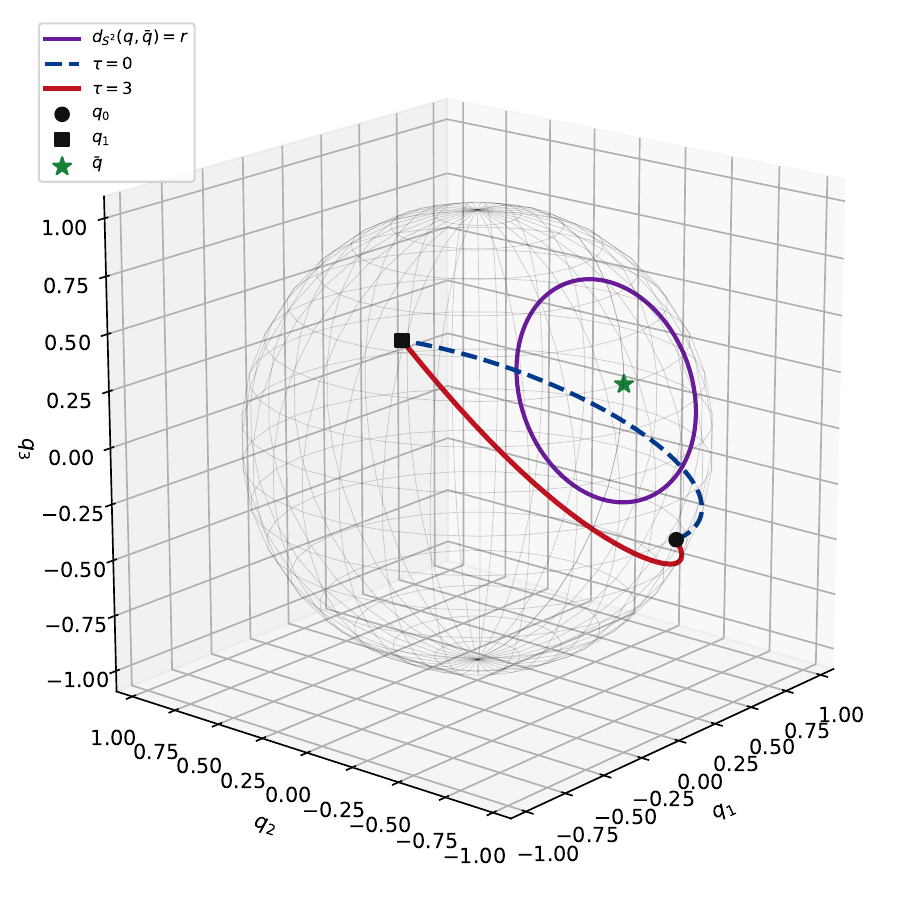} \caption{Projected camera-pointing trajectories. The shaded geodesic cap around the Sun direction $\bar q$ represents the undesirable pointing region. The dashed blue curve is the unpenalized baseline trajectory with $\tau=0$, which reaches a minimum distance of $8.992^\circ$ and enters the cap. The red curve is the penalized trajectory with $\tau=3$, whose minimum clearance increases to $36.031^\circ$.} \label{fig:satellite_projected_lgvi} 
\end{figure}

Figure~\ref{fig:satellite_projected_lgvi} shows the first comparison. The dashed blue curve is the baseline solution with $\tau=0$, using the same endpoints and the same quotient metric but no artificial potential penalty. The red curve is the solution with $\tau=3$. Without the artificial potential, the projected trajectory reaches a minimum distance of $8.992^\circ$ from $\bar q$ and therefore crosses the cap $\mathcal C_r(\bar q)$. With $\tau=3$, the minimum distance increases to $36.031^\circ$. Thus, although the formulation imposes only a soft penalty rather than a hard state constraint, the computed trajectory remains outside the undesirable pointing region. 

\begin{figure}[t] \centering 
\includegraphics[width=0.48\linewidth] {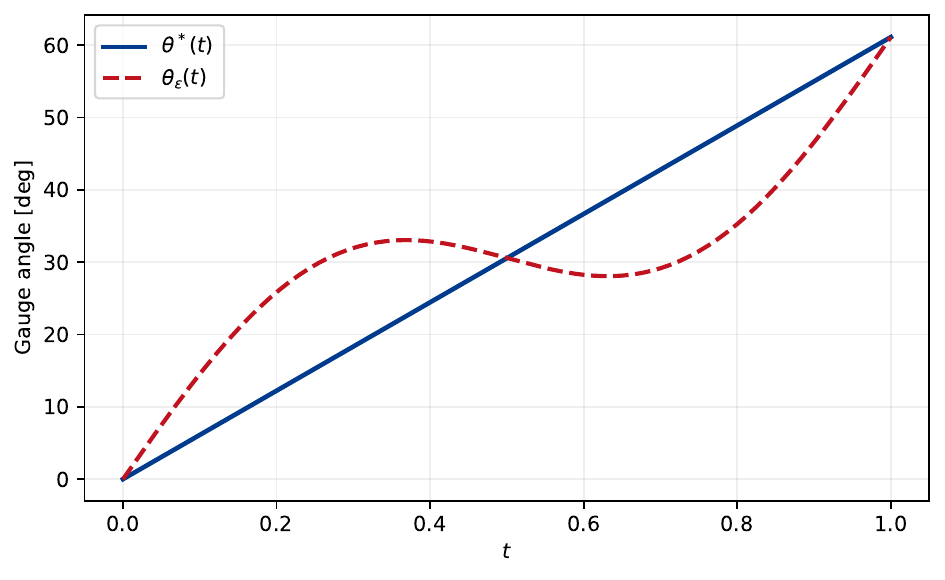} \hfill \includegraphics[width=0.48\linewidth] {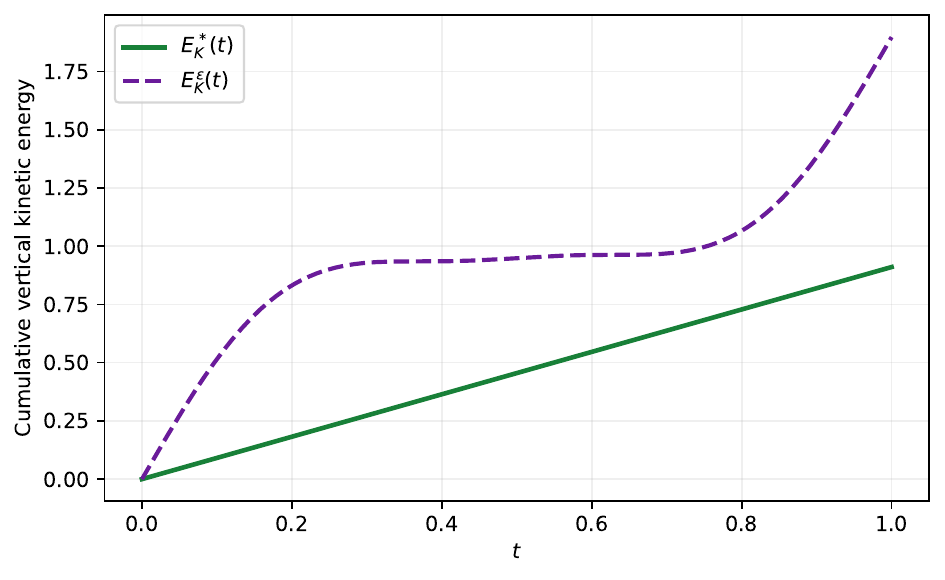} \caption{Optimal and comparison gauges with identical endpoint values and the same projected camera trajectory. Left: the optimal gauge $\theta^\ast(t)=\bar\theta t$ has constant angular rate, whereas $\theta_\varepsilon(t)$ introduces an additional oscillatory vertical motion. Right: the corresponding cumulative vertical kinetic energies. The constant-speed gauge yields $E_K[\kappa^\ast]=0.910397$, compared with $E_K[\kappa_\varepsilon]=1.89736$ for the perturbed gauge.} \label{fig:satellite_gauge_lgvi} 
\end{figure}

We now keep the penalized projected trajectory fixed and apply the optimal-gauge construction. For this trajectory, the terminal mismatch $\kappa_b=R_h(1)^{-1}R_1$ has principal angle $\bar\theta=61.121^\circ$, so the minimum-energy gauge is $\kappa^\ast(t)=\Exp(-t\bar\theta\,\hat e_3)$, or equivalently $\theta^\ast(t)=\bar\theta t$ in angular coordinates. To illustrate that the optimality is a statement about the vertical gauge, not about the projected trajectory, we compare this with the perturbed gauge $\theta_\varepsilon(t)=\bar\theta t+\varepsilon\sin(2\pi t)$, where $\varepsilon=0.25$. Both gauges have the same endpoint values and both lifts project to the same curve $q(t)$ on $S^2$, but the perturbed gauge introduces an unnecessary oscillatory rotation about the camera axis. Figure~\ref{fig:satellite_gauge_lgvi} shows this second comparison. The constant-speed gauge has vertical energy $E_K[\kappa^\ast]=0.910397$, whereas the oscillatory comparison gauge has $E_K[\kappa_\varepsilon]=1.89736$. This agrees with Theorem~\ref{thm: full_optimization_G}: among all lifts of the fixed projected trajectory satisfying the same attitude endpoints, the length-minimizing geodesic in $K\cong\SO(2)$ gives the minimum vertical kinetic energy. The discrete horizontal kinetic contribution is $E_{S^2}^d=3.23197$, and therefore the split kinetic energy of the optimal lift is $$E_{\SO(3)}^d=E_{S^2}^d+E_K[\kappa^\ast]=4.14237.$$ 
The LGVI preserves the group structure without reprojection. For $N=240$, the maximum error in $R_k^TR_k=I$ is $6.27\times10^{-15}$, the maximum determinant error is $5.00\times10^{-15}$, and the projected endpoint residual is $1.68\times10^{-15}$. The horizontal constraint is also preserved to numerical precision, with 
$$\max_k |(\mathbb J\Omega_k)_3|=1.71\times10^{-17}.$$
A refinement from $N=60$ to $N=480$ changes the minimum clearance from $36.051^\circ$ to $36.031^\circ$ and the optimal gauge angle from $61.124^\circ$ to $61.121^\circ$, which confirms numerical convergence of the reported quantities.

\section{Conclusions and Future Work}

We developed a variational lifting framework for mechanical systems on
Riemannian homogeneous spaces $H=G/K$. The original mechanical action on
$H$ is lifted to $G$ through a functional whose kinetic term depends only on
the horizontal velocity, and we showed that critical points of the two
variational problems are equivalent under projection. When the lifted
potential admits a partial symmetry, the corresponding dynamics admit a
reduced Euler--Poincar\'e description with an advected parameter. The
$H^1$ formulation makes explicit the degeneracy of the lifted problem:
although the projected critical curve and its horizontal representative are
smooth, an arbitrary $H^1$ curve in the isotropy subgroup may be appended
without changing the lifted action.

This gauge freedom led to a second variational problem selecting a
distinguished representative. Once the horizontal lift is fixed, minimizing
the total kinetic energy is equivalent to minimizing the vertical kinetic
energy on $K$. The optimal gauge is therefore a length-minimizing geodesic
in the isotropy subgroup, and the total kinetic energy splits into a projected
contribution and an explicit minimum vertical contribution. For closed
projected curves, the latter measures the minimum energy required to
compensate the holonomy of the horizontal lift. The construction was
illustrated on $\CP^n\cong SU(n+1)/S(U(1)\times U(n))$, where the gauge
freedom corresponds to the isotropy of unitary representatives of a pure
quantum state, and on $S^2\cong SO(3)/SO(2)$ through the camera-pointing
problem for an axisymmetric satellite. In the latter case, the numerical
results recover the projected avoidance trajectory and the constant-speed
minimum-energy gauge while preserving the Lie-group structure.

Several directions remain open. A natural next step is to develop a discrete
counterpart of the variational lifting and reduction procedure. The Lie group
variational integrator used in Section~4.4 preserves the $SO(3)$ configuration
space, but it discretizes the reduced rigid-body equations rather than the
general lifted variational principle developed here. It would therefore be
interesting to construct a discrete lifted action on $G$ whose invariance
under $K$-valued discrete gauges is preserved exactly, derive the associated
discrete reduced equations with symmetry breaking, and determine conditions
under which discretization and reduction commute. Such a construction should
also provide a discrete analogue of the horizontal--vertical energy splitting
and of the optimal-gauge problem on $K$.

Other directions include the treatment of endpoint conditions lying in
different connected components of $K$, the possible nonuniqueness of optimal
gauges at the cut locus, and the interaction between the intrinsic
minimum-energy gauge and restricted actuation. The latter is particularly
relevant for quantum and mechanical control systems, where the geometrically
optimal isotropy motion need not belong pointwise to the available control
distribution.

\section*{Acknowledgments}

The authors sincerely thank the anonymous referee for the careful reading of the manuscript and for the constructive comments that significantly improved both its presentation and mathematical clarity. 

L. J. Colombo was supported by Grant PID2022-137909NB-C21 funded by \\ MCIN/AEI/10.13039/501100011033. J.~R.~Goodman was supported by the Marie Sk\l odowska-Curie grant agreement No.~101206748 (GNACS).

\medskip

\textit{Use of artificial intelligence}. Generative AI tools were used to review and improve the writing of the manuscript and to assist with code for the numerical experiments. The authors take full responsibility for the content of the manuscript, including all mathematical statements, proofs, and numerical results.

\bibliographystyle{ws-gm}
\bibliography{references}

\label{eof}
\end{document}